\documentclass[11pt]{article}
\usepackage{amssymb}
\usepackage{latexsym}
\usepackage{amsfonts}
\usepackage{amsmath}
\newtheorem{thrm}{Theorem}[section]
\newtheorem{lemma}[thrm]{Lemma}
\newtheorem{prop}[thrm]{Proposition}

\newtheorem{cor}[thrm]{Corollary}
\newtheorem{remark}[thrm]{Remark}
\newtheorem{assumption}[thrm]{Assumption}
\newtheorem{example}[thrm]{Example}
\numberwithin{equation}{section}

\usepackage[dvips]{color}

\def\P{\mathbb{P} }
\def\R{\mathbb{R} }
\def\N{\mathbb{N} }
\def\V{\mathbb{V} }

\def\C{{\cal C} }

\begin{document}
\allowdisplaybreaks

\allowdisplaybreaks
\begin{doublespace}
\title{\Large\bf   Limit Theorems for Some Critical Superprocesses}
\author{ \bf  Yan-Xia Ren\footnote{The research of this author is supported by NSFC (Grant No.  11271030 and 11128101) and Specialized Research Fund for the
Doctoral Program of Higher Education.\hspace{1mm} } \hspace{1mm}\hspace{1mm}
Renming Song\thanks{Research supported in part by a grant from the Simons
Foundation (208236).} \hspace{1mm}\hspace{1mm} and \hspace{1mm}\hspace{1mm}
Rui Zhang\footnote{Supported by the China Scholarship Council. Corresponding author.}
\hspace{1mm} }
\date{}
\maketitle

\begin{abstract}
Let $X=\{X_t, t\ge 0; \P_\mu\}$ be a critical superprocess starting from a finite measure $\mu$.
Under some conditions, we first prove that
$\lim_{t\to\infty}t{\P}_{\mu}\left( \|X_t\|\ne 0 \right)=\nu^{-1}\langle \phi_0,\mu\rangle$,
where $\phi_0$ is the  eigenfunction corresponding to the first eigenvalue of the infinitesimal generator
$L$ of the mean semigroup of $X$, and $\nu$ is a positive constant.
Then we show that, for a large class of functions $f$,
conditioning on $\|X_t\|\ne 0$, $t^{-1}\langle f, X_t\rangle$
converges in distribution to $\langle f,\psi_0\rangle_m W$, where $W$
is an exponential random variable, and $\psi_0$ is the eigenfunction corresponding to
the first eigenvalue of  the dual of $L$. Finally, if $\langle f,\psi_0\rangle_m=0$,
we prove that, conditioning on $\|X_t\|\ne 0$, $\left( t^{-1}\langle \phi_0, X_t\rangle,
t^{-1/2}\langle f, X_t\rangle\right)$ converges in distribution
to $\left(W,G(f)\sqrt{W}\right)$, where $G(f)\sim\mathcal{N}(0,\sigma_f^2)$ is a normal random
variable, and $W$ and $G(f)$ are independent.
\end{abstract}


\medskip


\bigskip

\baselineskip=6.0mm
\section{Introduction}

\subsection{Motivation}

It is well known that if $\{Z_n, n\ge 0\}$ is a
critical (single type) branching process with finite second moment, then
\begin{equation}\label{extinction-descrete}
\lim_{n\to\infty}nP(Z_n>0)=\frac{2}{\sigma^2}
\end{equation}
and
\begin{equation}\label{clt-descrete}
\lim_{n\to\infty}P\left(\frac{1}{n}Z(n)>\frac{\sigma^2}{2}x|Z(n)>0\right)=e^{-x},\quad x\ge 0,
\end{equation}
where $\sigma^2$ is the variance of the offspring distribution.
The first result, due to Kolmogorov \cite{kol}, says that the non-extinction rate is of order $1/n$ as $n\to\infty$.
The second result, due to Yaglom \cite{yag}, says that
conditioned on non-extinction at time $n$,
the total population size in generation $n$ grows like  $n$.
For references to these results in English, one can see, for example,
\cite{harris} and \cite{KNS}.
For probabilistic proofs
of these results, see Lyons, Pemantle and Peres \cite{LPP}.
For continuous time critical branching processes $\{Z_t, t\ge 0\}$,
Athreya and Ney \cite[Theorem 3 and Lemma 2 on page 113]{AN72}
proved the following limit theorem:
\begin{equation}\label{clt-continuous}
\lim_{t\to\infty}P\left(\frac{1}{t}Z(t)>\frac{\sigma^2}{2}x|Z(t)>0\right)=e^{-x},\quad x\ge 0,
\end{equation}
where $\sigma^2$ is a positive constant determined by the branching rate and the variance
of the offspring distribution.

For discrete time multitype critical branching processes $\{{\bf Z}(n), n\ge 0\}$,
Athreya and Ney \cite{AN72}
gave two limit theorems under the finite second moment condition, see
\cite[Section V.5]{AN72}.
Let ${\bf v}$ be a positive left eigenvector of the mean matrix
associated with the eigenvalue 1.
The first order limit theorem says that
if ${\bf w}\cdot{\bf v}>0$, then
\begin{equation}\label{qusi-limit}
\lim_{n\to\infty}P\left(\frac{{\bf Z}(n)\cdot{\bf w}}{n}>x|{\bf Z}(n)>0\right)=
e^{-x/\gamma_1},
\quad x\ge 0,
\end{equation}
where $\gamma_1:=\gamma_1({\bf w})$ is a positive constant.
The second order limit theorem says that
if ${\bf w}\cdot{\bf v}=0$, then
\begin{equation}\label{central-limit}
\lim_{n\to\infty}P\left(\frac{{\bf Z}(n)\cdot{\bf w}}{\sqrt{n}}>x|{\bf Z}(n)>0\right)
=\int_x^\infty f(y)dy,
\quad x\in \R,
\end{equation}
where $$
f(y)=\frac{1}{2\gamma_2}e^{-|y|/\gamma_2},\quad y\in\R,
$$
and $\gamma_2:=\gamma_2({\bf w})$ is a positive constant.
The limit result \eqref{qusi-limit} is a generalization of \eqref{clt-descrete}
from the single type case to  the multitype case,
and was first proved by Joffe and Spitzer \cite{JS}.
The limit result \eqref{central-limit}
was first proved in Ney \cite{Ney}.

For continuous  time multitype critical branching processes,
Athreya and Ney \cite{AN74} proved two limit theorems,
similar to \eqref{qusi-limit} and \eqref{central-limit} respectively,
under the finite second moment condition,
see \cite[Theorems 1 and 2]{AN74}.

For limit theorems of critical branching processes (single type or multitype) without the finite second
moment condition, one can see, for instance \cite{GH, pak, sen, sla, vat, vat00} and the
references therein.

Asmussen and Hering \cite{AH83} discussed similar questions for critical  branching Markov processes
$\{Y_t,t\ge 0\}$ in a general space $E$ under the so-called condition (M)
(see \cite[page 156]{AH83}) on the first moment semigroup of $\{Y_t,t\ge 0\}$.
For each fixed $t\ge 0$, $Y_t$ is a random measure on $E$. For any finite measure $\mu$
on $E$ and any measurable function $f$ on $E$, we use $\|\mu\|$ to denote the total
mass of $\mu$ and $\langle f, \mu\rangle$ to denote the integral of $f$ with respect to $\mu$.
In \cite[Proposition 3.3 on page 201]{AH83}, Asmussen and Hering discussed the finite
time extinction property of branching Markov processes.
\cite[Theorem 3.4 on page 202]{AH83}
provided the rate of non-extinction, more precisely, it was shown that
 $$
\lim_{t\to\infty}tP_{\mu}(\|Y_t\|\neq 0)=\nu^{-1}\int_E\phi_0(x)\mu(dx)
 $$
uniformly in $\mu$ with $\|\mu\|=n$ for any integer positive $n$,
where $\nu$ is a positive constant and
$\phi_0$ is the first eigenfunction of the first moment  semigroup of $\{Y_t,t\ge 0\}$.
\cite[Theorem 3.8 on page 204]{AH83} gave a result similar to \eqref{qusi-limit},
while \cite[Theorem 3.3 on page 297]{AH83} gave a result similar to \eqref{central-limit}.
\cite[Section 4, Chapter VI]{AH83} discussed the limit theorems for the
critical branching Markov processes with infinite second moments.

As far as we know, not much has been done regarding limiting theorems of $\langle f, Y_t\rangle$
for critical branching Markov processes conditioned on $\|Y_t\|\neq 0$, since the book \cite{AH83}.
For critical superprocesses conditioned on non-extinction at time $t$,
Evans and Perkins \cite{Evans-Perkins} obtained results similar to
\eqref{extinction-descrete} and \eqref{qusi-limit} when $\varphi(x,z)=z^2$, $\beta(x)\equiv 1$
and the spatial process satisfies some ergodicity conditions.
\cite{Evans-Perkins} did not consider central limit theorem type results.
We note in passing that \cite{Evans-Perkins} also obtained results similar to
\eqref{qusi-limit}  conditioned on remote survival.
See \cite{Ch-Ro} for similar results for multitype Dawson-Watanabe processes
conditioned on remote survival.

The main purpose of this paper is to
establish limit theorems similar to \eqref{extinction-descrete},  \eqref{qusi-limit}
and a central limit type theorem
for critical superprocesses, under the finite second moment condition and
other very general, easy to check conditions.
Here is a summary of our main results.
Let $X=\{X_t, t\ge 0; \P_\mu\}$ be a critical superprocess starting from a finite measure $\mu$.
Under some conditions to be specified later, we first prove that
$\lim_{t\to\infty}t{\P}_{\mu}\left( \|X_t\|\ne 0 \right)=\nu^{-1}\langle \phi_0,\mu\rangle$,
where $\phi_0$ is the  eigenfunction corresponding to the first eigenvalue of the infinitesimal generator
$L$ of the mean semigroup of $X$, and $\nu$ is a positive constant.
Then we show that, for a large class of functions $f$,
conditioning on $\|X_t\|\ne 0$, $t^{-1}\langle f, X_t\rangle$
converges in distribution to $\langle f,\psi_0\rangle_m W$, where $W$
is an exponential random variable, and $\psi_0$ is the eigenfunction corresponding to
the first eigenvalue of  the dual of $L$. Finally, if $\langle f,\psi_0\rangle_m=0$,
we prove that, conditioning on $\|X_t\|\ne 0$, $\left( t^{-1}\langle \phi_0, X_t\rangle,
t^{-1/2}\langle f, X_t\rangle\right)$ converges in distribution
to $\left(W,G(f)\sqrt{W}\right)$, where $G(f)\sim\mathcal{N}(0,\sigma_f^2)$ is a normal random
variable, and $W$ and $G(f)$ are independent.

In our recent papers \cite{RSZ, RSZ2a},
we established some spatial centfral limit theorems for supercritical superprocesses.
See also \cite{AM, Mi, RSZ2, RSZ3} for related results for supercritical branching Markov processes and supercritical superprocesses.
Our original motivation for the present paper is to establish spatial central limit theorems for
critical superprocesses.
One of the main tools of the papers above is the analytical and spectral properties of
the Feynman-Kac semigroup of the spatial process,
which also play an important role in this paper.
We will assume that
the dual, with respect to a certain measure, of the semigroup of the spatial process
is a Markov semigroup. See the next subsection for details.

For branching Markov processes, there is a clear particle picture.
This particle structure was used essentially in proving the
central limit theorems for supercritical branching Markov processes
in \cite{AM, RSZ2, RSZ3}. For superprocesses, the particle picture
is less clear. In this case, the backbone decomposition or the excursion measures
are frequently used to describe the `infinitesimal particles'.
\cite{Mi,RSZ} used the backbone decomposition to establish central limit theorems for
supercritical super-OU processes, while \cite{RSZ2a} used the excursion measures of superprocesses
to prove central limit theorems for general supercritical superprocesses.
In this paper, we will also use the excursion measure to prove our central limit theorem.
Up to now, there is no known backbone decomposition for critical superprocesses conditioned on survival up to $t$ yet.

\subsection{Superprocesses and assumptions}

In this subsection,
we describe the superprocesses we are going to work with and
formulate our assumptions.
Since one of the main tools of this paper is the analytic properties
of the semigroup of the spatial process, we will need to assume that
the semigroup of the spatial process has a dual with respect to a certain
measure $m$ and the dual semigroup is Markovian.

Suppose that $E$ is a locally compact separable metric space and that $m$ is a $\sigma$-finite Borel measure
on $E$ with full support.
Suppose that $\partial$ is a separate point not contained
in $E$. $\partial$ will be interpreted as the cemetery point.
We will use $E_{\partial}$ to denote $E\cup\{\partial\}$.
Every function $f$ on $E$ is automatically extended to $E_{\partial}$ by setting $f(\partial)=0$.
 We will assume that $\xi=\{\xi_t,\Pi_x\}$ is a Hunt process on $E$ and $\zeta:=
 \inf\{t>0: \xi_t=\partial\}$ is the lifetime of $\xi$.
We will use $\{P_t:t\geq 0\}$ to denote the semigroup of $\xi$.
We will use
$\mathcal{B}(E)$ ($\mathcal{B}^+(E)$) to denote the
set of (non-negative)  Borel measurable functions on $E$, and use
$\mathcal{B}_b(E)$ ($\mathcal{B}_b^+(E)$) to denote the
set of (non-negative) bounded Borel measurable functions on $E$.

The superprocess $X=\{X_t:t\ge 0\}$ we are going to work with is determined by three parameters:
a spatial motion $\xi=\{\xi_t, \Pi_x\}$ on $E$
which is a Hunt process,
a branching rate function $\beta(x)$ on $E$ which is a non-negative bounded measurable function and a branching mechanism $\varphi$ of the form
\begin{equation}\label{e:branm}
\varphi(x,z)=-a(x)z+b(x)z^2+\int_{(0,+\infty)}(e^{-z y}-1+z y)n(x,dy),
\quad x\in E, \, z\ge 0,
\end{equation}
where $a\in \mathcal{B}_b(E)$, $b\in \mathcal{B}_b^+(E)$ and $n$ is a kernel from $E$ to $(0,\infty)$ satisfying
\begin{equation}\label{n:condition}
  \sup_{x\in E}\int_{(0,+\infty)} y^2 n(x,dy)<\infty.
\end{equation}
The assumption \eqref{n:condition} is the counterpart of the second moment condition in \cite{KNS}.
In the multitype continuous time branching process case, one does not need to take the supremum
and explicitly assume that $\beta$ is bounded, since the
state space of the spatial process, i.e., the type space, is finite.
Under this assumption, the superprocess $X$ has finite second moments (see \eqref{1.9} below).
In our paper, we will not consider the special case that
$\beta(\cdot)(b(\cdot)+n(\cdot,(0,\infty)))=0$, a.e.-m.

Let ${\cal M}_F(E)$
be the space of finite measures on $E$, equipped with topology of weak convergence.
The superprocess $X$ is a  Markov process taking values in ${\cal M}_F(E)$.
The existence of such superprocesses is well-known, see, for instance,
\cite{Dawson}, \cite{E.B.} or \cite{Li11}.
For any $\mu \in \mathcal{M}_F(E)$, we denote the
law
of $X$ with initial configuration $\mu$ by $\P_\mu$.
As usual,
$\langle f,\mu\rangle:=\int_E f(x)\mu(dx)$
and $\|\mu\|:=\langle 1,\mu\rangle$.
Throughout this paper, a real-valued function $u(t,x)$ on $[0,\infty)\times E_\partial$
is said to be locally bounded if for any $t>0$, $\sup_{s\in[0,t], x\in E_\partial}|u(s, x)|<\infty$.
According to \cite[Theorem 5.12]{Li11}, there is a Hunt process
$X=\{\Omega, {\cal G}, {\cal G}_t, X_t, \P_\mu\}$ taking values in  $\mathcal{M}_F(E)$,
such that, for every
$f\in \mathcal{B}^+_b(E)$ and $\mu \in \mathcal{M}_F(E)$,
\begin{equation}
  -\log \P_\mu\left(e^{-\langle f,X_t\rangle}\right)=\langle u_f(t,\cdot),\mu\rangle,
\end{equation}
where $u_f(t,x)$ is the unique locally bounded non-negative solution to the equation
\begin{equation}\label{u_f}
   u_f(t,x)+\Pi_x\int_0^t\Psi(\xi_s, u_f(t-s,\xi_s))ds=\Pi_x f(\xi_t), \quad x\in E_\partial,
\end{equation}
where $\Psi(x,z)=\beta(x)\varphi(x,z)$, $x\in E$ and $z\ge0$, while $\Psi(\partial,z)=0$, $z\ge0$.
Since $f(\partial)=0$, we have $u_f(t,\partial)=0$ for any $t\ge 0$.
In this paper, the superprocess we deal with is always this Hunt realization.
The function $\beta$ is usually called the branching rate, and $\varphi$ is
called the branching mechanism.
For more general superprocesses, $\beta$ can be a
measure on $E$. For the process $X$ in this paper, $\beta$ can be absorbed to $\varphi$.
That is to say, we could have, without loss of generality, supposed
that  $\beta\equiv 1$.
To be consistent with the formulations of our previous papers
\cite{RSZ, RSZ2, RSZ2a, RSZ3}, we keep $\beta$ as a function.

Define
\begin{equation}\label{e:alpha}
\alpha(x):=\beta(x)a(x)\quad \mbox{and}\quad
A(x):=\beta(x)\left( 2b(x)+\int_0^\infty y^2 n(x,dy)\right).
\end{equation}
Then, by our assumptions, $\alpha(x)\in \mathcal{B}_b(E)$ and
$A(x)\in\mathcal{B}_b^+(E)$.
Thus there exists $K>0$ such that
\begin{equation}\label{1.5}
  \sup_{x\in E}\left(|\alpha(x)|+A(x)\right)\le K.
\end{equation}
For any $f\in\mathcal{B}_b(E)$ and $(t, x)\in (0, \infty)\times E$, define
\begin{equation}\label{1.26}
   T_tf(x):=\Pi_x \left[e^{\int_0^t\alpha(\xi_s)\,ds}f(\xi_t)\right].
\end{equation}
It is well known that $T_tf(x)=\P_{\delta_x}\langle f,X_t\rangle$ for every $x\in E$.

Our standing assumption on $\xi$ is that there exists a family of
continuous strictly positive functions $\{p(t,x,y):t>0\}$ on $E\times E$ such that, for any $t>0$ and nonnegative function $f$ on $E$,
$$
  P_tf(x)=\int_E p(t,x,y)f(y)m(dy).
$$
Define
$$a_t(x):=\int_E p(t,x,y)^2\,m(dy),\qquad \hat{a}_t(x):=\int_E p(t,y,x)^2\,m(dy).$$
In this paper, we assume that
\begin{assumption}\label{assum1}
\begin{description}
  \item[(i)]  For any $t>0$, $\int_E p(t,x,y)\,m(dx)\le 1$.
  \item[(ii)] For any $t>0$, we have
  \begin{equation}\label{1.17}
    e_t:=\int_E a_t(x)\,m(dx)=\int_E\hat{a}_t(x)\,m(dx)=\int_E\int_E p(t,x,y)^2\,m(dy)\,m(dx)<\infty.
  \end{equation}
  Moreover, the functions $x\to a_t(x)$ and $x\to\hat{a}_t(x)$ are continuous on $E$ .
\end{description}
\end{assumption}

Note that, in Assumption \ref{assum1}(i), the integration is with respect to the first
space variable. It implies that the dual semigroup $\{{\widehat P}_t:t\ge 0\}$
of $\{P_t:t\ge 0\}$ with respect to $m$ defined by
$$
{\widehat P}_tf(x)=\int_Ep(t,y, x)f(y)m(dy)
$$
is Markovian. Assumption \ref{assum1}(ii) is a pretty weak $L^2$ condition and it allows
us to apply results on operator semigroups in Hilbert spaces.

By H\"older's inequality, we have
\begin{equation}\label{1.1}
 p(t+s,x,y)=\int_E p(t,x,z)p(s,z,y)\,m(dz)\le (a_t(x))^{1/2}(\hat{a}_s(y))^{1/2}.
\end{equation}

It is well known and easy to check that, $\{P_t:t\ge 0\}$
and $\{\widehat{P}_t:t\ge 0\}$ are strongly continuous
contraction semigroups on $L^2(E, m)$,
see \cite{RSZ3} for a proof.
Recall that $\{P_t:t\ge 0\}$ is a strongly continuous contraction semigroup
on $L^2(E, m)$ means that, for any $f\in L^2(E, m)$, $\lim_{t\to0}\|P_tf-f\|_2=0$
and $\|P_tf\|_2\le \|f\|_2$ for all $t\ge 0$.
We will use $\langle\cdot, \cdot\rangle_m$ to denote inner product in $L^2(E, m)$.
Since $p(t, x, y)$ is continuous in $(x, y)$,
by \eqref{1.1} and Assumption \ref{assum1}(ii),
using the dominated convergence theorem, we get that,
for any $f\in L^2(E,m)$,  $P_tf$  and $\widehat{P}_tf$ are continuous.

It follows from Assumption \ref{assum1}(ii) that,
for each $t>0$, $P_t$ and $\{\widehat{P}_t\}$ are compact operators on $L^2(E,m)$.
Let $\widetilde{L}$ and $\widehat{\widetilde{L}}$ be the infinitesimal generators of the semigroups $\{P_t\}$ and $\{\widehat{P}_t\}$ in $L^2(E,m)$ respectively.
Define $\widetilde{\lambda}_0:=\sup \Re(\sigma(\widetilde{L}))=\sup\Re(\sigma(\widehat{\widetilde{L}}))$.
By Jentzsch's theorem (Theorem V.6.6 on page 337 of \cite{Sch}),
$\widetilde{\lambda}_0$ is an eigenvalue of multiplicity 1 for both $\widetilde{L}$ and
$\widehat{\widetilde{L}}$, and that an eigenfunction $\widetilde{\phi}_0$ of $\widetilde{L}$ corresponding
to $\widetilde{\lambda}_0$ can be chosen to be strictly positive $m$-almost
everywhere with $\|\widetilde{\phi}_0\|_2=1$ and
an eigenfunction $\widetilde{\psi}_0$ of $\widehat{\widetilde{L}}$ corresponding to
$\widetilde{\lambda}_0$  can be chosen to be strictly positive $m$-almost
everywhere with $\langle \widetilde{\phi}_0, \widetilde{\psi}_0\rangle_m=1$.
Thus for $m$-almost every $x\in E$,
$$
e^{\widetilde{\lambda}_0}\widetilde{\phi}_0(x)=P_1\widetilde{\phi}_0(x),
\qquad
e^{\widetilde{\lambda}_0}\widetilde{\psi}_0(x)=\widehat{P}_1\widetilde{\psi}_0(x).
$$
Hence $\widetilde{\phi}_0$ and $\widetilde{\psi}_0$ can be chosen to be continuous and strictly
positive everywhere on $E$.

Our second assumption is

\begin{assumption}\label{assum3}
\begin{description}
 \item[(i)] $\widetilde{\phi}_0$ is bounded.

 \item[(ii)]
The semigroup $\{P_t,t\ge0\}$ is intrinsically ultracontractive,
that is, there exists $c_t>0$ such that
\begin{equation}\label{newcondition2}
  p(t,x,y)\le c_t\widetilde{\phi}_0(x)\widetilde{\psi}_0(y).
\end{equation}
\end{description}
\end{assumption}

Assumption \ref{assum3} is a pretty strong assumption on the semigroup
$\{P_t:t\ge0\}$. However, this assumption is satisfied in a lot of cases.
In Subsection \ref{ss:examples}, we will give many examples where
Assumptions \ref{assum1} and \ref{assum3} are satisfied. Here we only give one very
special example. If $E$ consists of finitely many points
and $\xi=\{\xi_t: t\ge 0\}$ is a conservative irreducible Markov process on $E$,
then $\xi$ satisfies Assumptions \ref{assum1} and \ref{assum3} for some finite measure $m$ on $E$ with full support.
So, as special cases, our results give the analogs of the results of  Athreya and Ney \cite{AN74}
for critical super-Markov chains.

We will prove in Lemma \ref{l:density} that
there exists a function $q(t, x, y)$ on $(0, \infty)\times E\times E$ which is continuous
in $(x, y)$ for each $t>0$ such that
\begin{equation}\label{comp0}
e^{-Kt}p(t,x,y) \le q(t,x,y)\le e^{Kt}p(t,x,y), \quad (t, x, y)\in (0, \infty)\times E\times E
\end{equation}
and that for any bounded Borel function $f$ and any $(t, x)\in (0, \infty)\times E$,
$$
  T_tf(x)=\int_E q(t,x,y)f(y)\,m(dy).
$$
It follows immediately that
\begin{equation}\label{Lp}
  \|T_tf\|_2\le e^{Kt}\|P_tf\|_2\le e^{Kt}\|f\|_2.
\end{equation}
In \cite{RSZ3}, we have proved that $\{T_t: t\ge 0\}$ is a strongly continuous semigroup on $L^2(E, m)$.
Let $\{\widehat{T}_t,t>0\}$ be the adjoint operators on $L^2(E,m)$ of $\{T_t,t>0\}$ , that is, for $f,g\in L^2(E,m)$,
$$\int_E f(x)T_tg(x)\,m(dx)=\int_E g(x)\widehat{T}_tf(x)\,m(dx)$$
and
$$\widehat{T}_tf(x)=\int_E q(t,y,x)f(y)\,m(dy).$$
We have proved in \cite{RSZ3} that $\{\widehat{T}_t: t\ge 0\}$ is also a strongly continuous semigroup on $L^2(E, m)$.
We claim that, for all $t>0$ and $f\in L^2(E,m)$, $T_tf$  and $\widehat{T}_tf$ are continuous.
In fact, since $q(t, x, y)$ is continuous in $(x, y)$, by \eqref{1.1}, \eqref{comp0} and Assumption 1.1(ii),
using the dominated convergence theorem, we get that, for any $f\in L^2(E,m)$,  $T_tf$  and $\widehat{T}_tf$ are continuous.

By Assumption \ref{assum1}(ii) and \eqref{comp0}, we get that
$$
\int_E\int_E q^2(t,x,y)\,m(x)\,m(dy)\le e^{2Kt}\int_E\int_E p^2(t,x,y)\,m(x)\,m(dy)<\infty.
$$
Thus, for each $t>0$, $T_t$ and $\{\widehat{T}_t\}$ are compact operators on $L^2(E,m)$.
Let $L$ and $\widehat{L}$ be the infinitesimal generators of the semigroups $\{T_t\}$ and $\{\widehat{T}_t\}$ in $L^2(E,m)$ respectively.
Define $\lambda_0:=\sup \Re(\sigma(L))=\sup\Re(\sigma(\widehat{L}))$.
By Jentzsch's theorem,
$\lambda_0$ is an eigenvalue of multiplicity 1 for both $L$ and
$\widehat{L}$, and that an eigenfunction $\phi_0$ of $L$ corresponding
to $\lambda_0$ can be chosen to be strictly positive $m$-almost
everywhere with $\|\phi_0\|_2=1$ and
an eigenfunction $\psi_0$ of $\widehat{L}$ corresponding to
$\lambda_0$  can be chosen to be strictly positive $m$-almost
everywhere with $\langle \phi_0,\psi_0\rangle_m=1$.
Thus for $m$-almost every $x\in E$,
$$
e^{\lambda_0}\phi_0(x)=T_1\phi_0(x),
\qquad
e^{\lambda_0}\psi_0(x)=\widehat{T}_1\psi_0(x).
$$
Hence $\psi_0$ and $\phi_0$ can be chosen to be continuous and strictly
positive everywhere on $E$.

Using Assumption \ref{assum3}, the boundedness of $\alpha$ and an
argument similar to that used in the proof of \cite[Theorem 3.4]{DS},
one can show the following:
\begin{description}
\item[(i)] $\phi_0$ is bounded.

\item[(ii)]
The semigroup $\{T_t,t\ge0\}$ is intrinsically ultracontractive,
that is, there exists $c_t>0$ such that
\begin{equation}\label{condition2}
  q(t,x,y)\le c_t\phi_0(x)\psi_0(y).
\end{equation}
\end{description}

The condition (M) on \cite[Page 156]{AH83} is a condition similar in spirit
to the intrinsic ultracontractivity of $\{T_t,t\ge0\}$. This condition is not very easy
to check. Essentially the only examples given in \cite{AH83} satisfying this
condition are branching diffusion processes in bounded smooth domains.
Our Assumption \ref{assum3} is in terms of
the intrinsic ultracontractivity of $\{P_t,t\ge0\}$.
Intrinsic ultracontractivity has been studied intensively in the last 30 years
and there are many results on the intrinsic ultracontractivity of semigroups.
Using these results, we will give in Subsection \ref{ss:examples} many examples
satisfying  Assumption \ref{assum3}.

Let $\lambda_\infty$ be the $L^\infty$-growth bound of the semigroup $T_t$, i.e.,
$$\lambda_\infty:=\lim_{t\to\infty}\frac{1}{t}\log\|T_t\|_{\infty,\infty}.$$ It is easy to see that $\lambda_0\le \lambda_\infty$. Note that $\lambda_0$ gives the rate of local growth when it is positive, and implies local
extinction otherwise. While if $\lambda_\infty\neq0$, then  in some sense,
the exponential growth/decay rate of $\langle 1,X_t\rangle$, the total mass of $X_t$, is $\lambda_\infty$, see \cite{ERS}.  According to \cite[Thorem 2.7]{KiSo08b}, under Assumptions \ref{assum1} and \ref{assum3},
there exist constants $\gamma>0$ and $c>0$ such that, for any $(t,x,y)\in (1,\infty)\times E\times E$, we have
\begin{equation}\label{density0}
  \left|e^{-\lambda_0t}q(t,x,y)-\phi_0(x)\psi_0(y)\right|\le ce^{-\gamma t}\phi_0(x)\psi_0(y).
\end{equation}
Hence for $(t,x,y)\in (1,\infty)\times E\times E$, we have
\begin{equation}\label{density2}e^{-\lambda_0t}q(t,x,y)\ge(1-ce^{-\gamma t})\phi_0(x)\psi_0(y).\end{equation}
Since $q(t,x,\cdot)\in L^1(E,m),$ we have
$\psi_0\in L^1(E,m)$.
Therefore, by \eqref{density0}, for $t>1$,
$\|T_t\|_{\infty,\infty}\le (1+c)\|\phi_0\|_\infty\langle 1, \psi_0\rangle_me^{\lambda_0t}$,
which implies $\lambda_0=\lambda_\infty$.

The main interest of this paper is on critical superprocesses, so we assume that
\begin{assumption}\label{assum2}
  $\lambda_0=0$.
\end{assumption}

Define $q_t(x):=\P_{\delta_x}(\|X_t\|=0)$.
Note that, since $\P_{\delta_x}\|X_t\|=T_t1(x)>0$, we have $\P_{\delta_x}(\|X_t\|=0)<1$.
In this paper, we also assume that
\begin{assumption}\label{assum4}
For any $t>0$ and $x\in E$,
  $q_t(x)\in(0,1).$
And, there exists $t_0>0$ such that,
\begin{equation}\label{4.2'}
  \inf_{x\in E} q_{t_0}(x)>0.
\end{equation}
\end{assumption}
In Subsection \ref{ss:extinc}, we will give
a sufficient condition (in term of the function $\Psi$) for Assumption \ref{assum4}.
In Lemma \ref{lem3.3}, we will show that, under our assumptions,
$\lim_{t\to\infty}q_t(x)=1,$ uniformly in $x\in E.$

\subsection{Main results}
In this subsection, we will state our main results.
In the following, we use the notation
$$\P_{t,\mu}(\cdot):=\P_{\mu}\left(\cdot\mid \|X_t\|\ne 0\right).$$
Recall that the process $X$ is defined on  $(\Omega,\mathcal{G})$.
Suppose that, for each $t>0$, $Y_t$ is a measurable map from  $(\Omega,\mathcal{G})$
to a Polish space $S$ and that
$Z$ is an $S$-valued random variable on a probability space $(\widetilde{\Omega}, \widetilde{\cal G}, P)$,
we write
$$
Y_t|_{\P_{t,\mu}}\stackrel{d}{\rightarrow}Z,
$$
if $\lim_{t\to\infty}\P_{t,\mu}[f(Y_t)]=P[f(Z)]$ for all bounded
continuous real-valued functions $f$ on $S$.

Define
\begin{equation}\label{nu}
  \nu:=\frac{1}{2}\langle A(\phi_0)^2,\psi_0\rangle_m.
\end{equation}
It is easy to see that $0<\nu<\infty$.
Define
$$
\mathcal{C}_p:=\{f\in\mathcal{B}(E): \langle|f|^p,\psi_0\rangle_m<\infty\}
$$
and $\C_p^+:=\C_p\cap \mathcal{B}^+(E)$.
By \eqref{condition2} and the fact that $q(t,x,y)$ is continuous,
using the dominated convergence theorem, we get that, for $f\in \C_1$, $T_tf(x)$ is continuous.
Since $\psi_0\in L^1(E,m)$, $\mathcal{B}_b(E)\subset\C_p$.
Moreover, by H\"{o}lder's inequality, we get $\C_2\subset\C_1$.

\begin{thrm}
\label{lem:2.1} For any non-zero $\mu\in\mathcal{M}_F(E)$,
\begin{equation}\label{2.4}
   \lim_{t\to\infty}t{\P}_{\mu}\left( \|X_t\|\ne 0 \right)=\nu^{-1}\langle \phi_0,\mu\rangle.
\end{equation}
Furthermore, the convergence above is uniform in $\mu$ with $\mu(E)\le M,$ where $M>0$ is any constant.
\end{thrm}

\begin{thrm}\label{The:2.2}
If $f\in\C_2$ then, for any non-zero $\mu\in\mathcal{M}_F(E)$, we have
\begin{equation}\label{the:1}
  t^{-1}\langle f,X_t\rangle |_{\P_{t,\mu}}\stackrel{d}{\rightarrow}\langle f,\psi_0\rangle_m W,
\end{equation}
where $W$ is an exponential random variable with parameter $1/\nu$.
In particular, we have
\begin{equation}\label{2.20}
  t^{-1}\langle \phi_0,X_t\rangle |_{\P_{t,\mu}}\stackrel{d}{\rightarrow} W.
\end{equation}
\end{thrm}

\begin{remark}
(1) The distributional limit $ \langle f,\psi_0\rangle_m W$ in Theorem \ref{The:2.2} does
not depend on the starting measure $\mu$.

(2)
Since $1\in \mathcal{B}_b(E)\subset\C_2$, thus
the limit result above implies that
$$
t^{-1}\langle 1,X_t\rangle |_{\P_{t,\mu}}\stackrel{d}{\rightarrow}\langle 1,\psi_0\rangle_m W,
$$
which says that, conditioned on no-extinction at time $t$, the growth rate of the total
mass $\langle 1,X_t\rangle$ is $t$ as $t\to\infty$.
\end{remark}

It is well known, see for instance \cite[Theorem A2.3]{Kall}, that the collection of Radon
measures on $E$
 equipped with the vague topology forms a Polish space.
Let $\rho(\cdot,\cdot)$ be a metric on the space of Radon measures on $E$ compatible with the vague topology.
Let $l$ be the finite (deterministic) measure on $E$ defined by $l(dx)=\psi_0(x)m(dx)$.

\begin{cor}\label{cor1.1}
For any $f\in\C_2$ and non-zero $\mu\in\mathcal{M}_F(E)$, it holds that, as $t\to\infty$,
\begin{equation}\label{3.1}
\frac{\langle f,X_t\rangle}{\langle \phi_0,X_t\rangle}|_{\P_{t,\mu}}\stackrel{d}{\rightarrow}
\langle f,\psi_0\rangle_m.
\end{equation}
Moreover, for any non-zero $\mu\in\mathcal{M}_F(E)$ and $\epsilon>0$,
$$
\lim_{t\to\infty}\P_{t, \mu}\left(\rho\left(\frac{X_t}{\langle \phi_0,X_t\rangle},l\right)\ge \epsilon\right)=0.
$$
\end{cor}

The above corollary can be thought of as a ``weak'' law of large numbers.
Thus it is natural to consider a corresponding central limit type theorem.
For this,  we need to find constants $a_t$ such that
$$
a_t\left(\frac{X_t}{\langle \phi_0, X_t\rangle}-l\right)|_{\P_{t,\mu}}
\stackrel{d}{\rightarrow} Y
$$
for some nontrivial finite random measure $Y$.
According to \cite[Theorem 16.16]{Kall}, it suffices to show that for each
continuous function $f$ with compact support in $E$,
$$
a_t\left(\frac{\langle f, X_t\rangle}{\langle \phi_0, X_t\rangle}-\langle f,\psi_0\rangle_m\right)
|_{\P_{t,\mu}}\stackrel{d}{\rightarrow} \langle f, Y\rangle.
$$
This is equivalent to finding $a_t$ such that
$$
a_t\frac{\langle \tilde f, X_t\rangle}{\langle \phi_0, X_t\rangle}
|_{\P_{t,\mu}}\stackrel{d}{\rightarrow} \langle f, Y\rangle,
$$
where $\tilde f=
f-\langle f, \psi_0\rangle_m \phi_0$ satisfies $\langle \tilde f, \psi_0\rangle_m=0$.
This is the reason that we consider only functions $f\in \C_2$ and $\langle f,\psi_0\rangle_m=0$
in the next theorem.

Define
\begin{equation}\label{sigma}
 \sigma_f^2=\int_0^\infty\langle A(T_s f)^2,\psi_0\rangle_m \,ds.
\end{equation}

\begin{thrm}\label{them2}
Suppose that $f\in \C_2$ and $\langle f,\psi_0\rangle_m=0$, then we have $\sigma_f^2<\infty$ and, for any non-zero $\mu\in\mathcal{M}_F(E)$,
\begin{equation}\label{3.9}
  \left(t^{-1}\langle\phi_0,X_t\rangle,t^{-1/2}\langle f,X_t\rangle\right)|_{\P_{t,\mu}}\stackrel{d}{\rightarrow}\left(W,G(f)\sqrt{W}\right),
\end{equation}
where $G(f)\sim\mathcal{N}(0,\sigma_f^2)$ is a normal random
variable and $W$ is the random variable defined in Theorem \ref{The:2.2}.
Moreover, $W$ and $G(f)$ are independent.
\end{thrm}

Combining Theorems \ref{The:2.2} and \ref{them2}, we see that,
when $\sigma_f^2>0$, the density of $G(f)\sqrt{W}$ is
$$
d(x)=\frac{1}{\sqrt{2\nu\sigma_f^2}}\exp\left\{-\frac{2|x|}{\sqrt{2\nu\sigma_f^2}}\right\}, \qquad x\in\R.
$$

As a consequence of Theorem \ref{them2}, we immediately get the following
result, which can be thought of as some sort of central limit theorem.

\begin{cor}\label{cor2}
Suppose that $f\in \C_2$ and $\langle f,\psi_0\rangle_m=0$, then we have $\sigma_f^2<\infty$ and, for any non-zero $\mu\in\mathcal{M}_F(E)$,
\begin{equation}\label{3.9'}
  \left(t^{-1}\langle\phi_0,X_t\rangle,\frac{\langle f,X_t\rangle}{\sqrt{\langle\phi_0,X_{t}\rangle}}\right)|_{\P_{t,\mu}}\stackrel{d}{\rightarrow}\left(W,G(f)\right),
\end{equation}
where $G(f)\sim\mathcal{N}(0,\sigma_f^2)$ is a normal random
variable and $W$ is the random variable defined in Theorem \ref{The:2.2}.
Moreover, $W$ and $G(f)$ are independent.
\end{cor}
\begin{remark}
{\rm
Suppose that $m$ is a probability measure, the spatial motion $\xi$ is conservative (that is, $P_t1=1$), and that
the branching mechanism is spatial-independent with
\begin{equation}\label{spatial-independent}
\Psi(z)=bz^2+\int_0^\infty(e^{-zy}-1+zy)n(dy),
\end{equation}
where $b\ge 0$ and $\int_0^\infty z^2 \,n(dz)<\infty$.
Then $T_t=P_t$,  $\lambda_0=0$ and $\phi_0(x)=1$.
Thus Assumption \ref{assum2} is satisfied.
The process $\{\|X_t\|,t\ge 0\}$ is a continuous state branching process with branching mechanism $\Psi$.
We assume that $\Psi$ satisfies the Grey condition:
\begin{equation}\label{Grey}
  \int^\infty \frac{1}{\Psi(z)}\,dz<\infty.
\end{equation}
Then, for any $\mu\in \mathcal{M}_F(E)$,
$$\lim_{t\to\infty} t\P_{\mu}(\|X_t\|\neq 0)=2A^{-1}\|\mu\|,$$
where $A=2b+\int_0^\infty y^2 n(dy)$,
and $$t^{-1}\|X_t\| |_{\P_{t,\mu}}\stackrel{d}{\rightarrow} W,$$
where $W$ is an exponential random variable with parameter $2A^{-1}$.
The proofs can be found in \cite{Lambert,Li00}.
It is easy to check that, under the assumptions above, Assumption \ref{assum4}
is satisfied, see the end of Subsection \ref{ss:extinc}.

Suppose that the spatial motion  $\xi$ satisfies Assumption \ref{assum1} and}

{\bf Assumption 1.2$'$}\quad There exists $t_0>0$ such that $a_{t_0}, \hat{a}_{t_0}\in L^2(E,m)$.

\noindent
{\rm Then using an argument similar to that in \cite[Lemma 2.6 (1)]{RSZ3},
we can get that, for $f\in L^2(E,m)\cap L^4(E,m)$,
$$\lim_{t\to\infty}\V ar_{\delta_x}\langle f,X_t\rangle=\sigma_f^2<\infty.$$
Thus, using the same arguments as in the proofs of Theorem \ref{The:2.2} and Theorem \ref{them2} below,
we can get that Theorem \ref{The:2.2} and Theorem \ref{them2} are also valid in this case
for $f\in L^2(E,m)\cap L^4(E,m)$ and $\mu\in \mathcal{M}_F(E)$ with compact support.
We will not give the detailed proof in this case.

Note that in this case we do not need Assumption \ref{assum3}.
One can check that super inward Ornstein-Uhlenbeck processes
satisfy Assumption \ref{assum1} and Assumption 1.2$'$,
see \cite[Examples 4.1]{CRSZ}.
Thus Theorem \ref{The:2.2} and Theorem \ref{them2} hold for
super inward Ornstein-Uhlenbeck processes with spatial-independent
branching mechanism $\Psi$ given by \eqref{spatial-independent}.}
\end{remark}

\subsection{Examples}\label{ss:examples}

In this subsection we present a list of examples which satisfy
Assumptions \ref{assum1} and \ref{assum3}.
For simplicity, we will not try to give the weakest possible conditions.
The first six are examples where the processes are symmetric with respect to some measure.

\begin{example} {\rm
Suppose that $E$ is a connected open subset of  $\mathbb{R}^d$ with finite Lebesgue measure
and that $m$ denotes the Lebesgue measure on $E$. Let $\xi$ be the subprocess in $E$ of a diffusion
process in $\mathbb{R}^d$ corresponding to a uniformly elliptic divergence form second order
differential operator. Then it is well known that $\xi$ has a transition density $p(t, x, y)$
which is a strictly positive, continuous and symmetric function of $(x, y)$ for any $t>0$
and that there exists $c>0$ such that
$$
p(t, x, y)\le c\,t^{d/2}, \qquad (t, x, y)\in (0, \infty)\times E\times E.
$$
Thus Assumption \ref{assum1} is trivially satisfied. If $E$ is a bounded Lipschitz connected
open set, then it follows from \cite{DS} that the semigroup $\{P_t:t\ge 0\}$ of $\xi$
is intrinsic ultracontractive and that the eigenfunction $\widetilde{\phi}_0$ corresponding
to the largest eigenvalue of the generator of $\{P_t:t\ge 0\}$ is bounded.
Thus Assumption \ref{assum3} is satisfied.
Under much weaker regularity
assumptions on $E$, Assumptions \ref{assum1} and \ref{assum3} are still satisfied. For some of these weaker
regularity assumptions, one can see \cite{ba} and the references therein.
}
\end{example}

\begin{example} {\rm
Suppose that $E$ is the closure of a bounded connected $C^2$ open set in $\mathbb{R}^d$ and that $m$
denotes the Lebesgue measure on $E$. Let $\xi$ be the reflecting Brownian motion
in $E$. Then $\xi$ has a transition density $p(t, x, y)$
which is a strictly positive, continuous and symmetric function of $(x, y)$ for any $t>0$
and that there exists $c>0$ such that
$$
p(t, x, y)\le c\,t^{d/2}, \qquad (t, x, y)\in (0, \infty)\times E\times E.
$$
The largest eigenvalue of the generator of the semigroup $\{P_t:t\ge 0\}$ of $\xi$
is $\widetilde{\lambda}_0=0$
and the corresponding eigenfunction $\widetilde{\phi}_0$ is a positive constant.
Thus Assumptions \ref{assum1} and \ref{assum3} are trivially satisfied.
}
\end{example}

\begin{example} {\rm
Suppose that $E$ is an open subset of  $\mathbb{R}^d$ with finite Lebesgue measure
and that $m$ denotes the Lebesgue measure on $E$.
Let $\xi$ be the subprocesses in $E$ of any of
the subordinate Brownian motions studied in \cite{KSV1, KSV2}.
Then it is known (see \cite{CKS2, CKS3}) that $\xi$ has a transition density $p(t, x, y)$
which is a strictly positive, continuous, bounded, symmetric function of $(x, y)$ for any $t>0$.
Thus Assumption \ref{assum1} is trivially satisfied. It follows from \cite{KiSo08c}
that the semigroup $\{P_t:t\ge 0\}$ of $\xi$
is intrinsic ultracontractive and that the eigenfunction $\widetilde{\phi}_0$ corresponding
to the largest eigenvalue of the generator of $\{P_t:t\ge 0\}$ is bounded.
Thus  Assumption \ref{assum3} is also satisfied.
}
\end{example}

\begin{example} {\rm
Suppose $a>2$ is a constant.
Assume that $E=\mathbb{R}^d$ and $m$ is the Lebesgue measure on $\mathbb{R}^d$.
Let $\xi$ be a Markov process on $\mathbb{R}^d$ corresponding to the infinitesimal
generator $\Delta-|x|^a$.
Let $p(t,x, y)$
denote the transition density of $\xi$ with respect to
the Lebesgue measure on $\mathbb{R}^d$.
It follows from \cite[Theorem 6.1]{DS} and its proof that, for any $t>0$, there exists $c_t>0$ such that
$$
p(t,x, y)\le c_t\exp\left(-\frac{2}{2+a}|x|^{1+a/2}\right)\exp\left(-\frac{2}{2+a}|y|^{1+a/2}\right),
\qquad x, y\in \mathbb{R}^d,
$$
that the eigenfunction $\widetilde{\phi}_0$ corresponding
to the largest eigenvalue of the generator of $\{P_t:t\ge 0\}$ of $\xi$
is bounded and that $\{P_t:t\ge 0\}$ is intrinsically ultracontractive.
Thus  Assumptions \ref{assum1} and \ref{assum3} are satisfied.
}
\end{example}

\begin{example}{\rm
Assume that $E=\mathbb{R}^d$ and $m$ is the Lebesgue measure on $\mathbb{R}^d$.
Suppose that $V$ is a nonnegative and locally bounded function on  $\mathbb{R}^d$
such that there exist $R>0$ and $M\ge 1$ such that for all $|x|>R$,
$$
M^{-1}(1+V(x))\le V(y)\le M(1+V(x)), \qquad y\in B(x, 1),
$$
and that
$$
\lim_{|x|\to\infty}\frac{V(x)}{\log|x|}=\infty.
$$
Suppose $\beta\in (0, 2)$ is a constant.
Let $\xi$ be a Markov process on $\mathbb{R}^d$ corresponding to the infinitesimal
generator $-(-\Delta)^{\beta/2}-V(x)$.
Let $p(t,x, y)$ denote the transition density of $\xi$ with respect to
the Lebesgue measure on $\mathbb{R}^d$.
It follows from \cite[Corollaries 3 and 4]{KK} that, for any $t>0$, there exists $c_t>0$ such that
$$
p(t,x, y)
\le c_t\frac1{(1+V(x))(1+|x|)^{d+\beta}}\frac1{(1+V(y))(1+|y|)^{d+\beta}},
\qquad x, y\in \mathbb{R}^d,
$$
that the eigenfunction $\widetilde{\phi}_0$ corresponding
to the largest eigenvalue of the generator of $\{P_t:t\ge 0\}$ of $\xi$
is bounded and that $\{P_t:t\ge 0\}$ is intrinsically ultraccontractive.
Thus  Assumptions \ref{assum1} and \ref{assum3} are satisfied.
}
\end{example}

\begin{example}{\rm
Assume that $E=\mathbb{R}^d$ and $m$ is the Lebesgue measure on $\mathbb{R}^d$.
A nondecreasing function $L:[0, \infty)\to [0, \infty)$ is said to be in the class ${\bf L}$
if $\lim_{t\to\infty}L(t)=\infty$ and there exists $c>1$ such that
$$
L(t+1)\le c(1+L(t)), \qquad t\ge 0.
$$
Suppose that $V$ is a nonnegative function on  $\mathbb{R}^d$
such that
$$
\lim_{|x|\to\infty}\frac{V(x)}{|x|}=\infty
$$
and that there exists a function $L\in {\bf L}$ such that there exists  $C>0$ such that
$$
L(|x|)\le V(x)\le C(1+L(|x|), \qquad x\in \mathbb{R}^d.
$$
Suppose that $r>0$
and $\beta\in (0, 2)$ are constants.
Let $\xi$ be a Markov process on $\mathbb{R}^d$ corresponding to the infinitesimal
generator $r- (-\Delta +r^{2/\beta})^{\beta/2}-V(x)$.
Let $p(t,x, y)$ denote the transition density of $\xi$ with respect to
the Lebesgue measure on $\mathbb{R}^d$.
It follows from \cite[Theorem 1.6]{KuSi} that, for any $t>0$, there exists $c_t>0$ such that
$$
p(t,x, y)\le c_t\frac{\exp(-r^{1/\beta}|x|)}{(1+V(x))(1+|x|)^{(d+\beta+1)/2}}
\frac{\exp(-r^{1/\beta}|y|)}{(1+V(y))(1+|y|)^{(d+\beta+1)/2}},
\qquad x, y\in \mathbb{R}^d,
$$
that the eigenfunction $\widetilde{\phi}_0$ corresponding
to the largest eigenvalue of the generator of $\{P_t:t\ge 0\}$ of $\xi$
is bounded and that $\{P_t:t\ge 0\}$ is intrinsically ultracontractive.
Thus  Assumptions \ref{assum1} and \ref{assum3} are satisfied.
}
\end{example}

In the next five examples the processes may not symmetric.

\begin{example}\label{examp1}
{\rm Suppose that $\beta\in (0, 2)$ and that $\xi^{(1)}=\{\xi^{(1)}_t: t\ge0\}$ is a strictly
$\beta$-stable process in $\mathbb{R}^d$. Suppose that, in the case $d\ge 2$, the spherical
part $\eta$ of the L\'evy measure $\mu$ of $\xi^{(1)}$ satisfies the following assumption: there
exist a positive function $\Phi$ on the unit sphere $S$ in $\mathbb{R}^d$ and $\kappa>1$
such that
$$
\Phi=\frac{d\eta}{d\sigma} \quad \mbox{and} \quad
\kappa^{-1}\le \Phi(z)\le \kappa \quad \mbox{on } S
$$
where $\sigma$ is the surface measure on $S$. In the case $d=1$, we assume that the L\'evy
measure of $\xi^{(1)}$ is given by
$$
\mu(dx)=c_1x^{-1-\beta}1_{\{x>0\}}+ c_2|x|^{-1-\beta}1_{\{x<0\}}
$$
with $c_1, c_2>0$. Suppose that $E$ is an open set in $\mathbb{R}^d$ of finite Lebesgue measure.
Let $\xi$ be the process in $E$ obtained by killing $\xi^{(1)}$ upon exiting $E$.
Then it follows from \cite[Example 4.1]{KiSo09} that $\xi$ has a transition density $p(t, x, y)$
which is a strictly positive, bounded continuous function of $(x, y)$ for any $t>0$.
Thus Assumption \ref{assum1} is trivially satisfied. It follows also from \cite[Example 4.1]{KiSo09}
that the semigroup $\{P_t:t\ge 0\}$ of $\xi$
is intrinsic ultracontractive and that the eigenfunction $\widetilde{\phi}_0$ corresponding
to the largest eigenvalue of the generator of $\{P_t:t\ge 0\}$ is bounded.
Thus Assumption \ref{assum3} is also satisfied.
}
\end{example}

\begin{example}\label{examp2}
{\rm Suppose that $\beta\in (0, 2)$ and that $\xi^{(2)}=\{\xi^{(2)}_t: t\ge 0\}$ is a truncated
strictly $\beta$-stable process in $\mathbb{R}^d$, that is, $\xi^{(2)}$ is a L\'evy
process with L\'evy measure given by
$$
\widetilde{\mu}(dx)=\mu(dx)1_{\{|x|<1\}},
$$
where $\mu$ is the L\'evy measure of the process $\xi^{(1)}$ in the previous example.
Suppose that $E$ is a connected open set in $\mathbb{R}^d$ of finite Lebesgue measure.
Let $\xi$ be the process in $E$ obtained by killing $\xi^{(2)}$ upon exiting $E$.
Then it follows from \cite[Example 4.2 and Proposition 4.4]{KiSo09}
that $\xi$ has a transition density $p(t, x, y)$
which is a strictly positive, bounded continuous function of $(x, y)$ for any $t>0$.
Thus Assumption \ref{assum1} is trivially satisfied. It follows also from
\cite[Example 4.2 and Proposition 4.4]{KiSo09}
that the semigroup $\{P_t:t\ge 0\}$ of $\xi$
is intrinsic ultracontractive and that the eigenfunction $\widetilde{\phi}_0$ corresponding
to the largest eigenvalue of the generator of $\{P_t:t\ge 0\}$ is bounded.
Thus Assumption \ref{assum3} is also satisfied.
}
\end{example}

\begin{example}\label{examp3}
{\rm
Suppose $\beta\in (0, 2)$, $\xi^{(1)}=\{\xi^{(1)}_t: t\ge0\}$ is a strictly
$\beta$-stable process in $\mathbb{R}^d$ satisfying the assumptions in Example
\ref{examp1} and that $B=\{B_t: t\ge 0\}$ is an independent Brownian motion in $\mathbb{R}^d$.
Let $\xi^{(3)}$ be the
process defined by $\xi^{(3)}_t=\xi^{(1)}_t+B_t$.
Suppose that $E$ is an open set in $\mathbb{R}^d$ of finite Lebesgue measure.
Let $\xi$ be the process in $E$ obtained by killing $\xi^{(3)}$ upon exiting $E$.
Then it follows from \cite[Example 4.5 and Lemma 4.6]{KiSo09}
that $\xi$ has a transition density $p(t, x, y)$
which is a strictly positive, bounded continuous function of $(x, y)$ for any $t>0$.
Thus Assumption \ref{assum1} is trivially satisfied. It follows also from
\cite[Example 4.5 and Lemma 4.6]{KiSo09}
that the semigroup $\{P_t:t\ge 0\}$ of $\xi$
is intrinsic ultracontractive and that the eigenfunction $\widetilde{\phi}_0$ corresponding
to the largest eigenvalue of the generator of $\{P_t:t\ge 0\}$ is bounded.
Thus Assumption \ref{assum3} is also satisfied.
}
\end{example}

\begin{example}\label{examp4}
{\rm
Suppose $\beta\in (0, 2)$, $\xi^{(2)}=\{\xi^{(2)}_t: t\ge0\}$ is a truncated strictly
$\beta$-stable process in $\mathbb{R}^d$ satisfying the assumptions in Example
\ref{examp2} and that $B=\{B_t: t\ge 0\}$ is an independent Brownian motion in $\mathbb{R}^d$.
Let $\xi^{(4)}$ be the
process defined by
$\xi^{(4)}_t=\xi^{(2)}_t+B_t$.
Suppose that $E$ is
a connected open set in $\mathbb{R}^d$ of finite Lebesgue measure.
Let $\xi$ be the process in $E$ obtained by killing $\xi^{(4)}$ upon exiting $E$.
Then it follows from \cite[Example 4.7 and Lemma 4.8]{KiSo09}
that $\xi$ has a transition density $p(t, x, y)$
which is a strictly positive, bounded continuous function of $(x, y)$ for any $t>0$.
Thus Assumption \ref{assum1} is trivially satisfied. It follows also from
\cite[Example 4.7 and Lemma 4.8]{KiSo09}
that the semigroup $\{P_t:t\ge 0\}$ of $\xi$
is intrinsic ultracontractive and that the eigenfunction $\widetilde{\phi}_0$ corresponding
to the largest eigenvalue of the generator of $\{P_t:t\ge 0\}$ is bounded.
Thus Assumption \ref{assum3} is also satisfied.
}
\end{example}

\begin{example}\label{examp5}
{\rm
Suppose $d\ge 3$ and that $\mu=(\mu^1, \cdots, \mu^d)$, where each $\mu^j$ is a signed measure
on
$\mathbb{R}^d$ such that
$$
 \lim_{r\to 0}\sup_{x\in\mathbb{R}^d}\int_{B(x, r)}\frac{|\mu^j|(dy)}{|x-y|^{d-1}}=0.
$$
Let $\xi^{(5)}=\{\xi^{(5)}_t: t\ge 0\}$ be a Brownian motion with drift $\mu$ in $\mathbb{R}^d$,
see \cite{KiSo06}.
Suppose that $E$ is a bounded connected open set in $\mathbb{R}^d$
and that $K>0$ is a constant such that
$E\subset B(0, K/2)$. Put $B=B(0, K)$. Let $G_B$ be the Green function of $\xi^{(5)}$ in $B$ and define
$H(x):=\int_BG_B(x, y)dy$. Then $H$ is a strictly positive continuous function on $B$.
Let $\xi$ be the
process obtained by killing $\xi^{(5)}$ upon exiting $E$.
Let $m$ be the measure on $E$ defined by $m(dx)=H(x)dx$.
Then it follows from \cite[Example 4.6]{ZLS} or \cite{KiSo08, KiSo08c}.
that $\xi$ has a transition density $p(t, x, y)$ with respect to $m$
and that $p(t, x, y)$ is a strictly positive, bounded continuous function of $(x, y)$ for any $t>0$.
Thus Assumption \ref{assum1} is trivially satisfied. It follows also from
\cite[Example 4.6]{ZLS} or \cite{KiSo08, KiSo08c}
that the semigroup $\{P_t:t\ge 0\}$ of $\xi$
is intrinsic ultracontractive and that the eigenfunction $\widetilde{\phi}_0$ corresponding
to the largest eigenvalue of the generator of $\{P_t:t\ge 0\}$ is bounded.
Thus Assumption \ref{assum3} is also satisfied.
}
\end{example}

\begin{example}\label{examp6}
{\rm
Suppose $d\ge 2$,
$\beta\in (1, 2)$,
and that $\mu=(\mu^1, \cdots, \mu^d)$, where each $\mu^j$ is a signed measure
on $\mathbb{R}^d$ such that
$$
 \lim_{r\to 0}\sup_{x\in\mathbb{R}^d}\int_{B(x, r)}\frac{|\mu^j|(dy)}{|x-y|^{d-\beta+1}}=0.
$$
Let $\xi^{(6)}=\{\xi^{(6)}_t: t\ge 0\}$ be an $\beta$-stable process with drift
$\mu$ in $\mathbb{R}^d$, see \cite{KiSo13}.
Suppose that $E$ is a bounded open set in $\mathbb{R}^d$ and suppose $K>0$ is such that
$D\subset B(0, K/2)$. Put $B=B(0, K)$. Let $G_B$ be the Green function of $\xi^{(6)}$ in $B$ and define
$H(x):=\int_BG_B(x, y)dy$. Then $H$ is a strictly positive continuous function on $B$.
Let $\xi$ be the
process obtained by killing $\xi^{(6)}$ upon exiting $D$.
Let $m$ be the measure on $E$ defined by $m(dx)=H(x)dx$.
Then it follows from \cite[Example 4.7]{ZLS} or \cite{CKS}
that $\xi$ has a transition density $p(t, x, y)$ with respect to $m$
and that $p(t, x, y)$ is a strictly positive,
bounded continuous function of $(x, y)$ for any $t>0$.
Thus Assumption \ref{assum1} is trivially satisfied. It follows also from
\cite[Example 4.7]{ZLS} or \cite{CKS}
that the semigroup $\{P_t:t\ge 0\}$ of $\xi$
is intrinsic ultracontractive and that the eigenfunction $\widetilde{\phi}_0$ corresponding
to the largest eigenvalue of the generator of $\{P_t:t\ge 0\}$ is bounded.
Thus Assumption \ref{assum3} is also satisfied.
}
\end{example}

\section{Preliminaries}

\subsection{Density of $\{T_t: t\ge 0\}$}

In this subsection, we show that, under Assumption \ref{assum1}, the semigroup $\{T_t: t\ge 0\}$
has a strictly positive density $q(t, x, y)$ and, for any $t>0$, $q(t, x, y)$ is continuous
in $(x, y)$.

\begin{lemma}\label{l:density}
Suppose that  Assumption \ref{assum1} holds. The semigroup $\{T_t: t\ge 0\}$ has a density $q(t, x, y)$
such that
\begin{equation}\label{comp}
 e^{-Kt}p(t,x,y) \le q(t,x,y)\le e^{Kt}p(t,x,y), \quad (t, x, y)\in (0, \infty)\times E\times E.
\end{equation}
Furthermore, for any $t>0$, $q(t, x, y)$ is a continuous function of $(x, y)$ on $E\times E$.
\end{lemma}

\textbf{Proof:}
For any $(t, x,y)\in (0, \infty)\times E\times E$, define
\begin{eqnarray*}
I_0(t, x, y)&:=&p(t,x, y),\\
I_n(t, x, y)&:=&\int^t_0\int_Ep(s,x, z)I_{n-1}(t-s, z, y)\alpha(z)m(dz)ds, \qquad n\ge 1.
\end{eqnarray*}
Using arguments similar to those in Subsection 1.2 of \cite{RSZ2}, we easily get that the function
\begin{equation}\label{e:denrel}
q(t,x, y):=\sum^\infty_{n=0}I_n(t, x, y), \qquad (t, x, y)\in (0, \infty)\times E\times E
\end{equation}
is well defined and $q(t,x,y)$ is the density of $T_t$ satisfying \eqref{comp}.
We omit the details.

We now prove the continuity of $q(t,x,y)$ in $(x,y)\in E\times E$ for each fixed $t>0$.
As in Subsection 1.2 of \cite{RSZ2}, it suffices to show that, for any $0<\epsilon<t/2$,
$$
\int^{t-\epsilon}_\epsilon \int_Ep(s,x, z)p(t-s,z, y)\alpha(z)m(dz)ds
$$
is continuous on $E\times E$.
By \eqref{1.1}, we get that
$$
p(s,x, z)p(t-s,z, y)|\alpha(z)| \le K
a_{\epsilon/2}(x)^{1/2}\widehat{a}_{\epsilon/2}(y)^{1/2}\widehat{a}_{s-\epsilon/2}(z)^{1/2}a_{t-s-\epsilon/2}(z)^{1/2}.
$$
We claim that the function $t\rightarrow \int_E a_{t}(x)\,m(dx)$ is decreasing.
Using this claim and H\"{o}lder's inequality, we get that
\begin{eqnarray*}
&&\int^{t-\epsilon}_\epsilon \int_E\widehat{a}_{s-\epsilon/2}(z)^{1/2}a_{t-s-\epsilon/2}(z)^{1/2}\,m(dz)\,ds\\
&\le& \int^{t-\epsilon}_\epsilon \left(\int_E\widehat{a}_{s-\epsilon/2}(z)\,m(dz)\right)^{1/2}\left(\int_Ea_{t-s-\epsilon/2}(z)\,m(dz)\right)^{1/2}\,ds \\
&=&
\int^{t-\epsilon}_\epsilon \left(\int_E a_{s-\epsilon/2}(z)\,m(dz)\right)^{1/2}\left(\int_Ea_{t-s-\epsilon/2}(z)\,m(dz)\right)^{1/2}\,ds
\le t\int_E a_{\epsilon/2}(z)\,m(dz).
\end{eqnarray*}
The  equality above follows from
the fact $\int_E \widehat{a}_{t}(z)\,m(dz)=\int_E a_{t}(z)\,m(dz)$.
Thus, by Assumption \ref{assum1}(ii) and the dominated convergence theorem, we get that
the function
$$
(x, y)\mapsto\int^{t-\epsilon}_\epsilon \int_Ep(s,x, z)p(t-s,z, y)\alpha(z)m(dz)ds
$$
is continuous.

Now, we prove the claim that the function $t\rightarrow \int_E a_{t}(x)\,m(dx)$ is decreasing.
In fact, by Fubini's theorem and  H\"{o}lder's inequality, we get
\begin{eqnarray*}
  a_{t+s}(x) &=& \int_E p(t+s,x,y)\int_{E}p(t,x,z)p(s,z,y)\,m(dz)\,m(dy)\\
&=& \int_E p(t,x,z)\int_{E}p(t+s,x,y)p(s,z,y)\,m(dy)\,m(dz)\\
&\le&  a_{t+s}(x)^{1/2}\int_E p(t,x,z)a_s(z)^{1/2}\,m(dz),
\end{eqnarray*}
which implies that
\begin{equation}\label{8.9}
  a_{t+s}(x)\le \left(\int_E p(t,x,z)a_s(z)^{1/2}\,m(dz)\right)^2\le \int_E p(t,x,z)a_s(z)\,m(dz).
\end{equation}
Thus, by Fubini's theorem and Assumption 1.1(i), we get that
\begin{equation}\label{8.10}
  \int_E a_{t+s}(x)\,m(dx)\le \int_E a_{s}(z)\int_E p(t,x,z)\,m(dx)\,m(dz)\le \int_E a_{s}(z)\,m(dz).
\end{equation}
We have now finished the proof of our claim.
\hfill$\Box$

\subsection{Extinction and non-extinction of $\{X_t,t\ge 0\}$}\label{ss:extinc}

In this subsection, we will give some sufficient
conditions for Assumption \ref{assum4},
see Lemma \ref{l:extinc} below. In the case when the function $a(x)$ in \eqref{e:branm}
is identically zero, this lemma
follows from \cite[Lemma 11.5.1]{Dawson}. Here
we provide a proof for completeness.

Let $\widetilde{\Psi}(x,z)$ be a function on $E_\partial\times(0,\infty)$ with the form:
 \begin{equation}\label{e:branm'}
\widetilde{\Psi}(x,z)=-\tilde{a}(x)z+\tilde{b}(x)z^2+\int_{(0,+\infty)}(e^{-z y}-1+z y)\tilde{n}(x,dy),
\quad x\in E_\partial, \quad z\ge 0,
\end{equation}
where $\tilde{a}\in \mathcal{B}_b(E_\partial)$, $\tilde{b}\in \mathcal{B}_b^+(E_\partial)$ and $\tilde{n}$ is a kernel from $E_\partial$ to $(0,\infty)$ satisfying
\begin{equation}\label{n:condition'}
   \int_{(0,+\infty)} (y\wedge y^2)\tilde{n}(x,dy)<\infty.
\end{equation}

The following Lemma \ref{l:comp} is similar to  \cite[Corollary 5.18 ]{Li11}.
Recall that, unless explicitly mentioned otherwise,
every function $f$ on $E$ is automatically extended to $E_{\partial}$ by setting $f(\partial)=0$.
The function $g$ in the lemma below may not satisfy $g(\partial)=0$.

\begin{lemma}\label{l:comp}
Suppose that $\Psi(x, z)\ge \widetilde{\Psi}(x, z)$ for all $x\in E$ and $z\ge0$.
Assume that $f$ and $g$
are bounded nonnegative measurable functions on $E_\partial$ such that $f(\partial)=0$ and
$f(x)\le g(x)$ for all $x\in E_\partial$.
If $v_g(t,x)$ is the unique locally bounded non-negative solution to the equation£º
$$
v_g(t,x)=-\Pi_x\int_0^t\widetilde{\Psi}(\xi_s, v_g(t-s, \xi_s))ds+\Pi_x g(\xi_t),\quad x\in E_{\partial}, t\ge0,
$$
then $v_g(t,x)\ge u_f(t,x)$ for all $t\ge0$ and $x\in E$, where $u_f$ is the unique
locally bounded non-negative solution to \eqref{u_f}.
\end{lemma}

\textbf{Proof:}
Recall that $u_f(t,\partial)=0$ and
$$
u_f(t,x)=-\int_0^t\Pi_x\big(\Psi(\xi_s, u_f(t-s, \xi_s)\big)ds+\Pi_x\big( f(\xi_t)\big),\quad x\in E_\partial.
$$
Define another branching mechanism $\Psi_1(x,z)$ as follows:
$$\Psi_1(x,z)=\left\{
                \begin{array}{ll}
                  \widetilde{\Psi}(x,z), & \hbox{$x\in E$;} \\
                  0, & \hbox{$x=\partial$.}
                \end{array}
              \right.
$$
Put $g_1(x)=g(x)1_E(x)$, for $x\in E_\partial$.
Then, for all $x\in E_\partial$, $\Psi_1(x,z)\le \Psi(x,z)$ and $f(x)\le g_1(x)$.
Let $u^1_{g_1}(t,x)$ be the unique locally bounded non-negative solution to the equation£º
$$
{u}^1_{g_1}(t,x)=-\int_0^t\Pi_x\big(\Psi_1(\xi_s, {u}^1_{g_1}(t-s, \xi_s))\big)ds+\Pi_x \big(g_1(\xi_t)\big),
\quad x\in E_{\partial}, t\ge0.
$$
It follows from \cite[Corollary 5.18 ]{Li11} that
\begin{equation}\label{ine1}u_f(t,x)\le {u}^1_{g_1}(t,x),\quad  x\in E,\, t\ge 0.\end{equation}
By \cite[Proposition 2.20 ]{Li11}, we have ${u}^1_{g_1}(t,\partial)\le
\Pi_{\partial} \left[e^{\int_0^t\alpha(\xi_s)\,ds}g_1(\xi_t)\right]=0$,
here we used the fact that $g_1(\partial)=0$. Therefore ${u}^1_{g_1}(t,\partial)=0$.
Since  $\widetilde{\Psi}(x,0)=0$, we have that
$$\Pi_x\big(\Psi_1(\xi_s, {u}^1_{g_1}(t-s, \xi_s))\big)=\Pi_x\big(\widetilde{\Psi}(\xi_s, {u}^1_{g_1}(t-s, \xi_s));s<\zeta\big)=\Pi_x\big(\widetilde{\Psi}(\xi_s, {u}^1_{g_1}(t-s, \xi_s))\big),$$
which implies that
${u}^1_{g_1}(t,x)$ is the unique locally bounded non-negative solution to the equation£º
$$
{u}^1_{g_1}(t,x)=-\int_0^t\Pi_x\big(\widetilde{\Psi}(\xi_s, {u}^1_{g_1}(t-s, \xi_s))\big)ds+\Pi_x \big(g_1(\xi_t)\big),\quad
x\in E_{\partial}, t\ge0.
$$
Since $g_1(x)\le g(x)$, for all $x\in E_\partial$, by \cite[Corollary 5.18]{Li11}, we have
\begin{equation}\label{ine2}
{u}^1_{g_1}(t,x)\le v_g(t,x)\quad  x\in E,\, t\ge 0.
\end{equation}
Combining \eqref{ine1} and \eqref{ine2}, we arrive at the desired assertion of this lemma.
\hfill$\Box$

\begin{lemma}\label{l:extinc}
Suppose that $\widetilde{\Psi}(z)\le\inf_{x\in E}\Psi(x,z)$, and $\widetilde{\Psi}(z)$ can be written in the form
$$
\widetilde{\Psi}(z)=\widetilde{a}z+\widetilde{b}z^2+\int^\infty_0
(e^{-zy}-1+zy)\widetilde{n}(dy)
$$
with $\widetilde{a}\in \R$, $\widetilde{b}\ge 0$ and $ \widetilde{n}$ is a measure on $(0,\infty)$ satisfying $\int^\infty_0(y\wedge y^2)\widetilde{n}(dy)<\infty$.
If $\widetilde{\Psi}(\infty)=\infty$
and $\widetilde{\Psi}(z)$ satisfies
\begin{equation}\label{Phi}
  \int^\infty\frac{1}{\widetilde{\Psi}(z)}\,dz<\infty,
\end{equation}
then, for any $t>0$,
$\inf_{x\in E}q_t(x)>0$.
\end{lemma}

\textbf{Proof:}
Let $\widetilde{X}$ be a continuous state branching
process with branching mechanism $\widetilde{\Psi}$.
 Let $\widetilde{\P}$ be the law of $\widetilde{X}$ with $\widetilde{X}_0=1$.
Define
$$
u_{\theta}(t, x)=-\log \P_{\delta_x} e^{-\theta \|X_t\|},
\qquad v_{\theta}(t)=-\log \widetilde{\P} e^{-\theta \widetilde{X}_t}.
$$
It is easy to see that $u_\theta(t,\partial)=0$ and,
for $x\in E$ and $t>0$,
$$
u_\theta(t,x)=-\Pi_x\int_0^t\Psi(\xi_s, u_f(t-s, \xi_s))ds+\theta\Pi_x (t<\zeta)
$$
and
$$
v_\theta(t)=-\int_0^t\widetilde{\Psi}(v_\theta(s))ds+\theta.
$$
Applying Lemma \ref{l:comp} with $\widetilde{\Psi}(x,z)=\widetilde{\Psi}(z)$, $x\in E_\partial,z\ge0$ and $g(x)=\theta,$ $x\in E_\partial$.
we get that, for all $t>0$, $x\in E$ and $\theta>0$,
$u_\theta(t,x)\le v_\theta(t)$.
Letting $\theta\to \infty$, we get
$-\log \mathbb{P}_{\delta_x}(\|X_t\|=0)\le -\log \widetilde{\mathbb{P}}(\tilde{X}_t=0)$.
 It is well known that, under the conditions of this lemma,
$\widetilde{\mathbb{P}}(\tilde{X}_t=0)>0$.
Thus
$\inf_{x\in E}q_t(x)=\inf_{x\in E}\mathbb{P}_{\delta_x}(\|X_t\|=0)\ge \widetilde{\mathbb{P}}(\tilde{X}_t=0)>0$.
\hfill$\Box$

\bigskip
It was proved in \cite{Sheu} that \eqref{Phi} is equivalent to $\int^\infty\frac{1}{\widetilde{\Phi}(z)}\,dz<\infty$,
where $\widetilde{\Phi}(z): =\widetilde{\Psi}(z)-\tilde a z$.
Lemma \ref{l:extinc} says that if the spatially dependent branching mechanism
$\Psi(x, z)$ is dominated from below by
a spatially independent branching mechanism $\widetilde\Psi(z)$ satisfying
$\displaystyle\widetilde\Psi(\infty)=\infty$ and \eqref{Phi}, then Assumption \ref{assum4} holds.
In particuler when $\Psi$
does not depend on the spatial variable $x$
and satisfies $\Psi(\infty)=\infty$ and the condition
$\int^\infty\frac{1}{\Psi(\lambda)}d\lambda<\infty$,
Assumption \ref{assum4} holds.
If $\tilde{b}:=\inf_{x\in E}b(x)\beta(x)>0$,
then $\Psi(x,z)\ge -K z+\tilde{b}z^2$, where $K$ is the constant given in \eqref{1.5}.
In this case, we can take $\widetilde{\Psi}(z):=-K z+\tilde{b}z^2$
and it is clear that $\widetilde{\Psi}(z)$ satisfies \eqref{Phi}.

\subsection{Estimates on moments}

In the remainder of this paper we will use the  following notation:
for two positive functions $f$ and $g$ on $E$, $f(x)\lesssim g(x)$ for $x\in E$ means that
there exists a constant $c>0$ such that
$f(x)\le cg(x)$ for all $x\in E$.
Throughout this paper, $c$ is a constant whose value may vary from line to line.

By \eqref{density0} and the assumption  that $\lambda_0=0$, we have,  for any $(t,x,y)\in (1,\infty)\times E\times E$,
\begin{equation}\label{density}
  \left|q(t,x,y)-\phi_0(x)\psi_0(y)\right|\le ce^{-\gamma t}\phi_0(x)\psi_0(y).
\end{equation}

It follows that, if $f\in \mathcal{C}_1$, we have, for $(t,x)\in(1,\infty)\times E$,
\begin{equation}\label{T_t}
  \left|T_tf(x)-\langle f,\psi_0\rangle_m \phi_0(x)\right|\le ce^{-\gamma t}\langle |f|,\psi_0\rangle_m \phi_0(x)
\end{equation}
 and
 \begin{equation}\label{1.6}
   \left|T_tf(x)\right|\le  (1+c)\langle |f|,\psi_0\rangle_m\phi_0(x).
 \end{equation}

Recall the second moment formula of the superprocess $\{X_t: t\ge 0\}$ (see, for example, \cite[Corollary 2.39]{Li11}):
for $f\in \mathcal{B}_b(E)$, we have for any $t >0$,
\begin{equation}\label{1.9}
  \P_{\mu}\langle f,X_t\rangle^2=\left(\P_{\mu}\langle f,X_t\rangle\right)^2+\int_E\int_0^tT_{s}[A(T_{t-s}f)^2](x)\,ds\mu(dx).
\end{equation}
Thus,
\begin{equation}\label{1.13}
   {\V}{\rm ar}_{\mu}\langle f,X_t\rangle=\langle{\V}{\rm ar}_{\delta_\cdot}\langle f,X_t\rangle, \mu\rangle=\int_E\int_0^tT_{s}[A(T_{t-s}f)^2](x)\,ds\mu(dx),
\end{equation}
where $\mathbb{V}{\rm ar}_{\mu}$ stands for the variance under $\P_{\mu}$.
For any $f\in \C_2$ and $x\in E$, applying the Cauchy-Schwarz inequality, we have $(T_{t-s}f)^2(x)\le e^{K(t-s)}T_{t-s}(f^2)(x)$, which implies that
\begin{equation}\label{4.4}
  \int_0^tT_{s}[A(T_{t-s}f)^2](x)\,ds\le e^{Kt}T_t(f^2)(x)<\infty.
\end{equation}
Thus, using a routine limit argument, one can easily check
that \eqref{1.9} and \eqref{1.13} also hold for $f\in \C_2$.

\begin{lemma} \label{lem:2.2}
Assume that $f\in \C_2$. If $\langle f,\psi_0\rangle_m=0$,
then for $(t, x)\in (2, \infty)\times E$, we have
\begin{equation}\label{2.10}
\left|{\V}{\rm ar}_{\delta_x}\langle f,X_t\rangle-\sigma^2_f\phi_0(x)\right|\lesssim e^{-\gamma t}\phi_0(x),
\end{equation}
where $\sigma_f^2$ is defined in \eqref{sigma}.
Therefore, for $(t, x)\in (2, \infty)\times E$, we have
\begin{equation}\label{2.11}
 {\V}{\rm ar}_{\delta_x}\langle f,X_t\rangle\lesssim \phi_0(x).
\end{equation}
\end{lemma}
\textbf{Proof:}
First, we show that $\sigma_f^2<\infty$.
For $s\le 1$, $|T_sf(x)|^2\le e^{Ks}T_s(f^2)(x)$.
Hence, for $s\le 1$,
\begin{equation}\label{6.1}
     \langle  A(T_s f)^2,\psi_0\rangle_m \le Ke^{sK}\langle  T_s(f^2),\psi_0\rangle_m=Ke^{sK}\langle f^2,\psi_0\rangle_m.
\end{equation}
For $s>1$, by \eqref{T_t}, $|T_sf(x)|\lesssim e^{-\gamma s}\langle |f|,\psi_0\rangle_m \phi_0(x)$.
Hence, for $s>1$,
\begin{equation}\label{6.2}
    \langle  A(T_s f)^2,\psi_0\rangle_m \lesssim e^{-2\gamma s}.
\end{equation}
Therefore, combining \eqref{6.1} and \eqref{6.2} we have that
$$
\sigma_f^2=\int_{0}^\infty \langle A(T_s f)^2,\psi_0\rangle_m \,ds\lesssim  \int_{0}^{1} e^{sK}\,ds+\int_{1}^\infty e^{-2\gamma s}\,ds<\infty.
$$
By \eqref{1.13}, for $t>2$, we have
\begin{eqnarray}
&&\left|\V{\rm ar}_{\delta_x}\langle f,X_t\rangle- \sigma^2_f \phi_0(x)\right|\nonumber\\
  &\le& \int^{t-1}_0\left|T_{t-s}[A(T_sf)^2](x)-\langle A(T_s f)^2,\psi_0\rangle_m\phi_0(x)\right|\,ds\nonumber\\
  &&+\int_{t-1}^tT_{t-s}[A(T_sf)^2](x)\,ds+\int_{t-1}^\infty \langle A(T_s f)^2,\psi_0\rangle_m \,ds \phi_0(x)\nonumber\\
  &=:& V_1(t,x)+V_2(t,x)+V_3(t,x).
\end{eqnarray}
First, we consider $V_1(t,x)$.
By \eqref{T_t}, for $t-s>1$, we have
\begin{eqnarray*}
  \left|T_{t-s}[A(T_sf)^2](x)-\langle A(T_s f)^2,\psi_0\rangle_m\phi_0(x)\right|
  &\lesssim &e^{-\gamma(t-s)}\langle A(T_s f)^2,\psi_0\rangle_m \phi_0(x).
\end{eqnarray*}
Therefore, by \eqref{6.1} and \eqref{6.2}, we have, for $(t,x)\in(2,\infty)\times E$,
\begin{eqnarray}\label{V1}
  V_1(t,x) &\lesssim&
  \int_{1}^{t} e^{-\gamma(t+s)}\,ds\,\phi_0(x)
  +\int_{0}^{1}e^{-\gamma(t-s)}\,ds\,\phi_0(x)
   \lesssim e^{-\gamma t}\phi_0(x).
\end{eqnarray}
For $V_2(t,x)$, by \eqref{T_t}, for $s>t-1>1$, $|T_sf(x)|\lesssim e^{-\gamma s}\phi_0(x)$.
Thus,
\begin{equation}\label{1.14}
  V_2(t,x)\lesssim \int_{t-1}^t e^{-2\gamma s}T_{t-s}[\phi_0^2](x)\,ds=e^{-2\gamma t}\int_{0}^{1} e^{2\gamma s}T_{s}[\phi_0^2](x)\,ds.
\end{equation}
By  H\"{o}lder's inequality, we have
$$
\phi_0^2(x)=(T_{1}\phi_0(x))^2\le e^{K}T_{1}(\phi_0^2)(x).$$
Thus by \eqref{1.14} and \eqref{1.6}, for $(t,x)\in(2,\infty)\times E$, we have
\begin{equation}\label{V_2}
  V_2(t,x)\lesssim e^{-2\gamma t}\int_{0}^{1} T_{s+1}(\phi_0^2)(x)\,ds\lesssim e^{-2\gamma t}\phi_0(x).
\end{equation}
For $V_3(t,x)$, by \eqref{T_t}, for $s>t-1>1$, $|T_sf(x)|\lesssim e^{-\gamma s}\phi_0(x)$.
Thus,
\begin{eqnarray}\label{1.15}
  V_3(t,x) &\lesssim& \int_{t-1}^\infty e^{-2\gamma s}\,ds\langle \phi_0^2,\psi_0\rangle_m \phi_0(x)
   \lesssim   e^{-2\gamma t}\phi_0(x).
\end{eqnarray}
It follows from \eqref{V1}, \eqref{V_2} and \eqref{1.15} that, for $(t,x)\in(2,\infty)\times E$,
$$
\left|\V{\rm ar}_{\delta_x}\langle f,X_t\rangle- \sigma^2_f \phi_0(x)\right|\lesssim e^{-\gamma t}\phi_0(x).
$$
Now \eqref{2.11} follows immediately.
\hfill$\Box$

\subsection{Excursion measures of $\{X_t,t\ge 0\}$}

We use $\mathbb{D}$ to denote the space of $\mathcal{M}_F({E})$-valued
right continuous functions $t\mapsto \omega_t$ on $(0, \infty)$ having zero as a trap.
We use $(\mathcal{A},\mathcal{A}_t)$ to denote the natural $\sigma$-algebras on
$\mathbb{D}$ generated by the coordinate process.

Under Assumption \ref{assum4}, it is known (see \cite[Chapter 8]{Li11}) that
one can associate with $\{\P_{\delta_x}:x\in E\}$ a family
of $\sigma$-finite measures $\{\N_x:x\in E\}$ defined on $(\mathbb{D},\mathcal{A})$ such  that $\mathbb{N}_x(\{0\})=0$, and
\begin{equation*}
\mathbb{N}_x (1- e^{-\langle f, \omega_{t} \rangle})
= -\log \mathbb{P}_{\delta_x}(e^{-\langle f, X_{t} \rangle}) ,
\quad f\in {\cal B}^+_b(E),\ t\ge 0.
\end{equation*}
For further information on excursion measures of superprocesses,
we refer the reader to \cite{E.B2., elk-roe, Li03}.

For any $\mu\in {\cal M}_F(E)$, let $N(d\omega)$
be a Poisson random measure on the space $\mathbb{D}$
with intensity $\int_ E\mathbb{N}_x(d\omega)\mu(dx)$, in a probability space $(\widetilde\Omega, \widetilde{\cal F}, \mathbf{P}_\mu)$.
We define another process $\{\Lambda_t: t\ge 0\}$ by $\Lambda_0=\mu$ and
$$\Lambda_t:=\int_{\mathbb{D}}\omega_tN(d\omega),\quad t>0.$$
Let $\widetilde{\mathcal{F}}_t$ be the $\sigma$-algebra generated by the random variables $\{N(A):A\in \mathcal{A}_t\}$.
Then $\{\Lambda,(\widetilde{\mathcal{F}}_t)_{t\ge 0}, \mathbf{P}_{\mu}\}$ has the same law as
$\{X,(\mathcal{G}_t)_{t\ge 0},\P_{\mu}\}$, see \cite[Theorem 8.24]{Li11}.

Now we list some properties of $\mathbb{N}_x$. The proofs are similar to those of \cite[Corollary 1.2, Proposition 1.1]{E.B2.}.
\begin{prop}\label{Moments of N}
If $\mathbb{P}_{\delta_x}|\langle f, X_t\rangle|<\infty,$ then
\begin{equation}\label{N1}
 \mathbb{N}_x \langle f, \omega_t\rangle=\mathbb{P}_{\delta_x}\langle f, X_t\rangle.
\end{equation}
If $\mathbb{P}_{\delta_x}\langle f, X_t\rangle^2<\infty,$ then
\begin{equation}\label{N2}
 \mathbb{N}_x \langle f, \omega_t\rangle^2=\mathbb{V}ar_{\delta_x}\langle f, X_t\rangle.
\end{equation}
\end{prop}

For $f\in \mathcal{C}_1$, by \eqref{1.6}, $T_tf$ is bounded and in $\mathcal{C}_1$.
It follows from Proposition \ref{Moments of N} that, for any $f\in \mathcal{C}_1$,
$$
\int_E\int_D\langle |f|, \omega_s\rangle\mathbb{N}_x(d\omega)X_t(dx)<\infty,
\quad  \mathbb{P}_\mu\mbox{-a.s.}
$$
Now, by the Markov property of $X$, we get that for any $f\in \mathcal{C}_1$,
\begin{eqnarray}\label{cf}
  \mathbb{P}_\mu\left[\exp\left\{i\theta\langle f, X_{t+s}\rangle\right\}|X_t\right]
 & =&\mathbb{P}_{X_t}\left[\exp\left\{i\theta\langle f, X_{s}\rangle\right\}\right]=\textbf{P}_{X_t}\left[\exp\left\{i\theta\langle f, \Lambda_{s}\rangle\right\}\right]\nonumber\\
  &=&\exp\left\{\int_{E}\int_{\mathbb{D}}(e^{i\theta\langle f, \omega_s\rangle}-1)\mathbb{N}_x(d\omega)X_t(dx)\right\}.
\end{eqnarray}

\section{Proofs of Main Results}
In this section, we will prove our main theorems.

\subsection{Proofs of Theorems \ref{lem:2.1} and \ref{The:2.2}}
For $x\in E$ and $z>0$, define
\begin{equation}\label{def-r}r(x,z)=\Psi(x,z)+\alpha(x)z\end{equation} and
\begin{equation}\label{def-r2}r^{(2)}(x,z)=\Psi(x,z)+\alpha(x)z-\frac{1}{2}A(x)z^2.\end{equation}

\begin{lemma}\label{lem:1.3}  For any $x\in E$ and $z>0$,
\begin{equation}\label{2.6}
  0\le r(x,z)\le Kz^2/2
\end{equation}
and
\begin{equation}\label{2.14}
  |r^{(2)}(x,z)|\le e(x,z)z^2,
\end{equation}
where
\begin{equation}\label{de:e}
  e(x,z)=\beta(x)\int_0^\infty y^2\left(1\wedge\frac{1}{6}yz\right)\,n(x,dy).
\end{equation}
\end{lemma}
\textbf{Proof:}
It is easy to see that
\begin{equation}\label{2.2}
  r(x,z)=\beta(x)\left(b(x)z^2+\int_0^\infty(e^{-zy}-1+zy)\,n(x,dy)\right)
\end{equation}
and
$$r^{(2)}(x,z)=\beta(x)\int_0^\infty\left(e^{-zy}-1+zy-\frac{1}{2}y^2z^2\right)\,n(x,dy).$$
It follows from Taylor's expansion that, for $\theta>0$,
\begin{equation}\label{2.7}
  0<e^{-\theta}-1+\theta\le \frac{1}{2}\theta^2
\end{equation}
and
\begin{equation}\label{2.8}
  \left|e^{-\theta}-1+\theta- \frac{1}{2}\theta^2\right|\le\frac{1}{6}\theta^3 .
\end{equation}
By \eqref{2.7}, we also have $\left|e^{-\theta}-1+\theta- \frac{1}{2}\theta^2\right|\le \theta^2$.
Thus, we have
\begin{equation}\label{2.12}
  \left|e^{-\theta}-1+\theta- \frac{1}{2}\theta^2\right|\le \theta^2\left(1\wedge\frac{1}{6}\theta\right).
\end{equation}
Therefore, by \eqref{2.7} and \eqref{2.12}, we have
$$0<r(x,z)\le \beta(x)\left(b(x)+\frac{1}{2}\int_0^\infty y^2\,n(x,dy)\right)z^2\le Kz^2/2$$
and
$$r^{(2)}(x,z)\le\beta(x)\int_0^\infty y^2\left(1\wedge\frac{1}{6}yz\right)\,n(x,dy)z^2.$$
The proof is now complete.
\hfill$\Box$

\bigskip
Recall that $$u_f(t,x):=-\log\P_{\delta_x}e^{-\langle f,X_t\rangle}.$$

\begin{lemma}\label{lem1.2}
If $f\in \C_1^+$,
then $0\le u_f(t,x)<\infty$ for all $t\ge 0, x \in E$, and the function $R_f$ defined by
\begin{eqnarray}
 R_f(t,x)&:=&T_tf(x)-u_f(t,x)\label{de:R1}
\end{eqnarray}
satisfies
\begin{equation}\label{IF}
R_f(t,x)=\int_0^t T_s\left[r(\cdot,u_f(t-s,\cdot)\right](x)\,ds,\quad t\ge 0, x\in E.
\end{equation}
Moreover,
\begin{equation}\label{1.31}
  0\le R_f(t,x)\le e^{Kt}T_t(f^2)(x),\quad t\ge 0, x\in E.
\end{equation}
\end{lemma}
\textbf{Proof:}
First, we assume that $f\in \mathcal{B}_b^+$.
Recall that
$u_f(t,x)=-\log\P_{\delta_x}e^{-\langle f,X_t\rangle}$ satisfies
$$
  u_f(t,x)+\Pi_x\int_0^t \Psi(\xi_s,u_f(t-s,\xi_s))\,ds=\Pi_x(f(\xi_t)),\quad t\ge 0, x\in E.
$$
It follows from \cite[Theorem 2.23]{Li11} that $u_f(t, x)$ also satisfies
\begin{equation}\label{1.8}
  u_f(t,x)=-\int_0^t T_s\left[r(\cdot,u_f(t-s,\cdot))\right](x)\,ds+T_tf(x), \quad t\ge 0, x\in E.
\end{equation}
Thus, we get \eqref{IF} immediately.

For general $f\in \C_1^+$,  we have $T_tf(x)<\infty$.
Let $f_n(x)=f(x)\wedge n\in \mathcal{B}_b^+$.
Since \eqref{IF} holds for $f_n$,
applying the monotone convergence theorem,
we get that \eqref{IF} also holds for $f$.
Therefore, by \eqref{2.6}, $R_f(t,x)\ge 0$, which means $u_f(t,x)\le T_tf(x)<\infty$.
Recall that, as a consequence of the Cauchy-Schwarz inequality, we have
$(T_{t-s}f)^2(y)\le e^{K(t-s)}T_{t-s}(f^2)(y)$.
Combining this with \eqref{2.6}, we get
$$
  0\le R_f(t,x)\le \frac{K}{2}\int_0^t T_s[(u_f(t-s))^2](x)\,ds\le \frac{K}{2}\int_0^t T_s[(T_{t-s}f)^2](x)ds
  \le e^{Kt}T_t(f^2)(x).
$$
\hfill$\Box$

\bigskip
Recall that $q_t(x)=\P_{\delta_x}(\|X_t\|=0)$.
\begin{lemma}\label{lem3.3}
\begin{equation}\label{4.2}
   \lim_{t\to\infty}\inf_{x\in E}q_t(x)=1.
\end{equation}
\end{lemma}
\textbf{Proof:}
For $\theta>0$, let
$$
  u_\theta(t,x):=-\log \P_{\delta_x}e^{-\langle \theta,X_t\rangle}.
$$
By the Markov property of $X$,
\begin{equation}\label{1.20}
  q_{t+s}(x)=\lim_{\theta\to\infty}\P_{\delta_x}\left(e^{-\theta\|X_{t+s}\|}\right)
  =\lim_{\theta\to\infty}\P_{\delta_x}\left(e^{-\langle u_\theta(s)),X_t\rangle}\right)
  =\P_{\delta_x}\left(e^{-\langle -\log q_s ,X_t\rangle}\right).
\end{equation}
Since $q_t(x)$ is increasing in $t$, $q(x):=\lim_{t\to\infty}q_t(x)$ exists.
Put $w(x)=-\log q(x)$.
Letting $s\to \infty$ in \eqref{1.20},
we get $q(x)=\P_{\delta_x}\left(e^{-\langle w ,X_t\rangle}\right),$
which implies, for $t>0$,
\begin{equation}\label{new:1}
  w(x)=u_{w}(t,x)\quad x\in E.
\end{equation}
By Assumption 1.4,
for $s>t_0$,
$$\|w\|_\infty\le \|-\log q_{s}\|_\infty\le \|-\log q_{t_0}\|_\infty=-\log\left(\inf_{x\in E}q_{t_0}(x)\right)<\infty,$$
which implies $w\in \C_1^+$, and $ -\log q_s\in\C_1^+$.
Thus, by \eqref{de:R1}, \eqref{IF} and \eqref{new:1}, we have
\begin{equation}\label{1.12}
  w(x)=T_t(w)(x)-\int_0^t T_s(r(\cdot,w(\cdot)))(x)\,ds,\quad x\in E.
\end{equation}
By \eqref{T_t}, we have $\lim_{t\to\infty}T_t(w)(x)=\langle w,\psi_0\rangle_m\phi_0(x).$

If $\langle r(\cdot, w(\cdot)),\psi_0\rangle_m>0$, then
$$\lim_{t\to\infty}T_t(r(\cdot,w))(x)=\langle r(\cdot, w(\cdot)),\psi_0\rangle_m\phi_0(x)>0,\quad \mbox{for any }x\in E,$$
which implies
$$\lim_{t\to\infty}\int_0^t T_s\left[r(\cdot,w(\cdot))\right](x)\,ds=\infty, \quad \mbox{for any }x\in E.$$
Thus, by \eqref{1.12}, we get
$$0\le w(x)=\lim_{t\to\infty}\left(T_t(w))(x)-\int_0^t T_s\left[r(\cdot,w(\cdot))\right](x)\,ds\right)=-\infty,$$
which is a contradiction.
Therefore $r(x,w(x))=0$, a.e.-m.
Then, by \eqref{1.12}, we get, for all $x\in E$,
\begin{equation}\label{1.40}
  w(x)=\langle w,\psi_0\rangle_m\phi_0(x),
\end{equation}
which implies that $w\equiv 0$ on $E$ or $w(x)>0$ for any $x\in E$.
Since $r(x,w(x))=0$, a.e.-m., by \eqref{2.2}, we obtain  $w\equiv0$ on $E$.
For $s>t_0$, by \eqref{1.20} and Lemma \ref{lem1.2}, we get
$$-\log q_{2+s}(x)=u_{-\log q_s}(2,x)\le T_2(-\log q_s)(x)\le (1+c)\langle-\log q_s,\psi_0\rangle_m\|\phi_0\|_\infty,$$
where in the last inequality we used \eqref{T_t}.
Since $-\log q_s(x)\to 0$, by the dominated convergence theorem, we get
$$\lim_{s\to\infty}\langle -\log q_s,\psi_0\rangle_m=0.$$
Now \eqref{4.2} follows immediately.\hfill$\Box$

\begin{lemma}\label{lem:2.3}
For any $f\in\C_1^+$, there exists a function $h_f(t,x)$ such that
\begin{equation}\label{1.10}
  u_f(t,x)=(1+h_f(t,x))\langle u_f(t,\cdot),\psi_0\rangle_m\phi_0(x).
\end{equation}
Furthermore,
\begin{equation}\label{1.22}
  \lim_{t\to\infty}\|h_f(t)\|_\infty=0 \quad \mbox{uniformly in }  f\in \C_1^+.
\end{equation}
\end{lemma}
\textbf{Proof:}
For any $f\in \C_1^+$, we have $u_f(t,x)\le T_t f(x)< \infty$ and $\langle u_f(t,\cdot), \psi_0\rangle_m\le \langle T_tf, \psi_0\rangle_m=\langle f,\psi_0\rangle_m<\infty$. So $u_f(t,x)\in \C_1^+$.
If $m(f>0)=0$,
then $T_tf(x)=0$ for all $t>0$ and $x\in E$, which implies
$u_f(t,x)=0$ and $\langle u_f(t,\cdot),\psi_0\rangle_m=0$.
In this case, we define $h_f(t,x)=0$.
If $m(f>0)>0$, then $T_tf(x)>0$ for all $t>0$ and $x\in E$,
which implies $\P_{\delta_x}\left(\langle f,X_t\rangle=0\right)<1$.
Thus we have $u_f(t,x)>0$ and $\langle u_f(t,\cdot),\psi_0\rangle_m>0$.
Define
$$h_f(t,x)=\frac{ u_f(t,x)-\langle u_f(t,\cdot), \psi_0\rangle_m\phi_0(x)}{\langle u_f(t,\cdot), \psi_0\rangle_m\phi_0(x)}.$$
We only need to prove that  $\|h_f(t,\cdot)\|_\infty\to 0$ uniformly in $f\in \C_1^+\setminus\{0\}$ as $t\to\infty$.
Since $\P_{\mu}\big(e^{-\langle f,X_t\rangle}\big)\ge \P_{\mu}(\|X_t\|=0)$, we get that
\begin{equation}\label{2.5}
  \|u_f(t,\cdot)\|_\infty\le \|-\log q_t\|_\infty\to 0\quad \mbox{as } t\to\infty.
\end{equation}
By the Markov property of $X$ we have
\begin{eqnarray}\label{omega-prog}
  u_f(t,x)=-\log \P_{\delta_x}e^{-\langle u_f(t-s,\cdot),X_s\rangle}=u_{u_f(t-s)}(s, x),\quad t\ge s>0, x\in E.
\end{eqnarray}
where in the subscript on the right-hand side, $u_f(t-s)$ stands for the function $x\to u_f(t-s, x)$. In the remainder
of this proof, we keep this convention.
By \eqref{de:R1}, we have
\begin{equation}
 u_f(t,x)=T_s(u_f(t-s,\cdot))(x)-
R_{u_f(t-s)}(s, x).
 \label{1.21}
\end{equation}
Thus,
\begin{equation}\label{1.30}
  \langle u_f(t,\cdot),\psi_0\rangle_m=\langle u_f(t-s,\cdot),\psi_0\rangle_m-
   \langle R_{u_f(t-s)}(s, \cdot),\psi_0\rangle_m.
\end{equation}
Therefore, by \eqref{T_t}, \eqref{1.6} and \eqref{1.31}, we have, for $1<s<t$ and $x\in E$,
\begin{eqnarray*}
  &&|u_f(t,x)-\langle u_f(t,\cdot), \psi_0\rangle_m\phi_0(x)|\\
  &&\le|T_s(u_f(t-s,\cdot))(x)-\langle u_f(t-s,\cdot),\psi_0\rangle_m\phi_0(x)|\\
 &&\qquad+\left|R_{u_f(t-s)}(s, x)\right|+
 \left|\langle R_{u_f(t-s)}(s, \cdot),\psi_0\rangle_m\phi_0(x)\right|\\
 &&\le  ce^{-\gamma s}\langle u_f(t-s,\cdot),\psi_0\rangle_m\phi_0(x)
 +e^{Ks}T_s(u^2_f(t-s,\cdot))(x)+e^{Ks}\langle u^2_f(t-s,\cdot),\psi_0\rangle_m\phi_0(x)\\
  &&\le  ce^{-\gamma s}\langle  u_f(t-s,\cdot),\psi_0\rangle_m \phi_0(x)
 +(2+c)e^{Ks}\langle u^2_f(t-s,\cdot),\psi_0\rangle_m\phi_0(x)\\
  &&\le  \left[ce^{-\gamma s}+(2+c)e^{Ks}\|-\log q_{t-s}\|_\infty\right]\langle u_f(t-s,\cdot),\psi_0\rangle_m \phi_0(x),
\end{eqnarray*}
where in the last inequality we used \eqref{2.5} and $c$ is the constant in \eqref{density}.

By Lemma \ref{lem1.2} and \eqref{1.21}, we get
\begin{eqnarray}
 &&T_s(u_f(t-s,\cdot))(x)\ge
  u_f(t,x)\ge T_s(u_f(t-s,\cdot))(x)-e^{Ks}T_s(u^2_f(t-s,\cdot))(x)\nonumber\\
 &&\ge T_s(u_f(t-s,\cdot))(x)-e^{Ks}\|-\log q_{t-s}\|_\infty T_s(u_f(t-s,\cdot))(x).\label{1.11}
\end{eqnarray}
Thus,  we have
\begin{equation}\label{1.23}
  \langle u_f(t-s,\cdot),\psi_0\rangle_m\ge
  \langle u_f(t,\cdot), \psi_0\rangle_m \ge (1-e^{Ks}\|-\log q_{t-s}\|_\infty)\langle u_f(t-s,\cdot),\psi_0\rangle_m.
\end{equation}
For any $s>1$, $(1-e^{Ks}\|-\log q_{t-s}\|_\infty)>0$ when $t$ is large enough. Therefore, as $t\to\infty$,
$$
 \|h_f(t,\cdot)\|_\infty
 \le \frac{ce^{-\gamma s}+(2+c)e^{Ks}\|-\log q_{t-s}\|_\infty}{1-e^{Ks}\|-\log q_{t-s}\|_\infty}\to ce^{-\gamma s}.
$$
Now, letting $s\to\infty$, we get $\|h_f(t,\cdot)\|_\infty\to 0$ uniformly in $f\in \C_1^+\setminus\{0\}$ as $t\to\infty$.
\hfill$\Box$

\bigskip

\begin{lemma}\label{le:1.4}
For any $\delta>0$,
\begin{equation}\label{2.3}
  \lim_{n\to\infty}\frac{1}{n\delta}\left(\frac{1}{\langle u_{f}(n\delta,\cdot),\psi_0\rangle_m}-\frac{1}{\langle f,\psi_0\rangle_m}\right)=\nu
\end{equation}
uniformly in $f\in\C_1^+\setminus\{0\}$. Here $\nu$ is defined in \eqref{nu}.
\end{lemma}
\textbf{Proof:}
In this proof, we sometimes use $u_f(t)$ to denote the function $x\to u_f(t,x)$.
Since $f$ is non-negative and $m(f>0)>0$, we have
$u_f(t,x)>0$ for all $t>0$ and $x\in E$. Consequently, we have
$\langle u_f(t),\psi_0\rangle_m>0$.
It is clear that $u_f(0)=f$. First note that
\begin{eqnarray*}
  &&\frac{1}{n\delta}\left(\frac{1}{\langle u_f(n\delta),\psi_0\rangle_m}-\frac{1}{\langle f,\psi_0\rangle_m}\right) \\
   &=& \frac{1}{n\delta}\sum_{k=0}^{n-1}\left(\frac{1}{\langle u_f((k+1)\delta),\psi_0\rangle_m}-\frac{1}{\langle u_f(k\delta),\psi_0\rangle_m}\right) \\
   &=&\frac{1}{n\delta}\sum_{k=0}^{n-1}\left(\frac{\langle u_f(k\delta),\psi_0\rangle_m-\langle u_f((k+1)\delta),\psi_0\rangle_m}
   {\langle u_f((k+1)\delta),\psi_0\rangle_m\langle u_f(k\delta),\psi_0\rangle_m}\right).
\end{eqnarray*}
Recall the identity \eqref{omega-prog} and the definition of $r^{(2)}(x,z)$
given in \eqref{def-r2}. Using \eqref{1.30} with $t=(k+1)\delta$ and $s=\delta$, we get
\begin{eqnarray*}
 &&\langle u_f(k\delta), \psi_0\rangle_m-\langle u_f((k+1)\delta),\psi_0\rangle_m
 =\langle R_{u_f(k\delta)}(\delta, \cdot),\psi_0\rangle_m\\
   &&= \int_0^{\delta}\langle r(\cdot,u_f(k\delta+s,\cdot)),\psi_0\rangle_m\,ds\\
  &&=\frac{1}{2}\int_0^{\delta} \langle A(u_f(k\delta+s))^2,\psi_0\rangle_m\,ds
  +\int_0^{\delta}\langle r^{(2)}(\cdot,u_f(k\delta+s,\cdot)),\psi_0\rangle_m\,ds\\
  &&=:I_1+I_2.
\end{eqnarray*}
By \eqref{1.10} and \eqref{1.23}, we have, for $s\in[0,\delta]$,
\begin{eqnarray}\label{2.26}
  &&|u_f(t+s,x)-\langle u_f(t),\psi_0\rangle_m\phi_0(x)|\nonumber\\
  &\le& |u_f(t+s,x)-\langle u_f(t+s),\psi_0\rangle_m\phi_0(x)|
  +|\langle u_f(t),\psi_0\rangle_m-\langle u_f(t+s),\psi_0\rangle_m|\phi_0(x)\nonumber \\
  &\le& \|h_{f}(t+s)\|_\infty \langle u_f(t+s),\psi_0\rangle_m\phi_0(x)+e^{Ks}\|-\log q_{t}\|_\infty\langle u_f(t),\psi_0\rangle_m\phi_0(x)\nonumber \\
  &\le&\left(\|h_{f}(t+s)\|_\infty
  +e^{Ks}\|-\log q_{t}\|_\infty\right)\langle u_f(t),\psi_0\rangle_m\phi_0(x)\nonumber \\
 &\le &c_f(t)\langle u_f(t),\psi_0\rangle_m\phi_0(x),
\end{eqnarray}
where $c_f(t)=\sup_{0\le s\le \delta}\left(\|h_{f}(t+s)\|_\infty+e^{Ks}\|-\log q_{t}\|_\infty\right).$
By \eqref{1.22} and Lemma \ref{lem3.3}, we get
$c_f(t)\to 0$, as $t\to\infty$, uniformly in $f\in\C_1^+$.
Thus, by \eqref{2.26} we have for $s\in[0,\delta]$,
\begin{eqnarray}
 &&\frac{|u_f(t+s,x)^2-\langle u_f(t),\psi_0\rangle_m^2(\phi_0(x))^2|}{\langle u_f(t),\psi_0\rangle_m^2} \le\left(2+c_f(t)\right)c_f(t)(\phi_0(x))^2.\label{1.29}
\end{eqnarray}
Therefore, we have,
\begin{eqnarray*}
  \left|\frac{I_1}{\langle u_f(k\delta),\psi_0\rangle_m^2} -\delta\nu\right|
  &=& \frac{\left|\int_0^{\delta} \langle A\left((u_f(k\delta+s))^2-\langle u_f(k\delta),\psi_0\rangle_m^2\phi_0^2\right),\psi_0\rangle_m\,ds\right|}
  {2\langle u_f(k\delta),\psi_0\rangle_m^2} \\
  &\le&\frac{1}{2}\langle A\phi_0^2,\psi_0\rangle_m\delta \left(2+c_f(k\delta)\right)c_f(k\delta) \to 0,\quad \mbox{as}\quad k\to \infty,
\end{eqnarray*}
uniformly in $f\in \C_1^+\setminus\{0\}$.
By \eqref{1.23}, we have
$$
  0\le 1-\frac{\langle u_f((k+1)\delta),\psi_0\rangle_m}{\langle u_f(k\delta),\psi_0\rangle_m}\le e^{K\delta}\|-\log q_{k\delta}\|_\infty,
$$
which implies that
\begin{equation}\label{1.27}
  \frac{\langle u_f(k\delta),\psi_0\rangle_m}{\langle u_f((k+1)\delta),\psi_0\rangle_m}\to 1,
  \quad \mbox{ as } k\to \infty,
\end{equation}
uniformly in $f\in\C_1^+\setminus\{0\}$.
It follows that
\begin{equation}\label{1.28}
  \lim_{k\to\infty}\frac{I_1}{\langle u_f(k\delta),\psi_0\rangle_m \langle u_f((k+1)\delta),\psi_0\rangle_m} =\delta\nu
\end{equation}
uniformly in $f\in\C_1^+\setminus\{0\}$.

For $I_2$, by \eqref{2.14} and \eqref{2.26}, we have
\begin{eqnarray*}
 \frac{\langle r^{(2)}(\cdot,u_f(k\delta+s,\cdot)),\psi_0\rangle_m}{\langle u_f(k\delta),\psi_0\rangle_m^2}
  &\le& \frac{\langle e(\cdot,u_f(k\delta+s,\cdot))u_f(k\delta+s)^2,\psi_0\rangle_m}{\langle u_f(k\delta),\psi_0\rangle_m^2}\\
  &\le&  (1+c_f(k\delta))^2\langle e(\cdot,u_f(k\delta+s,\cdot))\phi_0^2,\psi_0\rangle_m \\
   &\le& (1+c_f(k\delta))^2\langle e(\cdot,\|-\log q_{k\delta}\|_\infty)\phi_0^2,\psi_0\rangle_m,
\end{eqnarray*}
here the last inequality follows from $\|u_f(k\delta+u)\|_\infty\le \|-\log q_{k\delta+u}\|_\infty\le \|-\log q_{k\delta}\|_\infty$ and the fact $z\rightarrow e(x,z)$ is increasing.
It is easy to see that the function $e(x,z)\downarrow 0$ as $z\downarrow0$.
Thus, as $k\to\infty$,
$$\frac{I_2}{\langle u_f(k\delta),\psi_0\rangle_m^2}
\le \delta (1+c_f(k\delta))^2\langle e(\cdot,\|-\log q_{k\delta}\|_\infty)\phi_0^2,\psi_0\rangle_m\to 0$$
uniformly in $f\in \C_1^+\setminus\{0\}$.
By \eqref{1.27}, we have
\begin{equation}\label{2.30}
  \lim_{k\to\infty}\frac{I_2}{\langle u_f(k\delta),\psi_0\rangle_m \langle u_f((k+1)\delta),\psi_0\rangle_m}=0
\end{equation}
uniformly in $f\in \C_1^+\setminus\{0\}$.
Using \eqref{1.28} and \eqref{2.30}, we get,
$$
\lim_{k\to\infty}\frac{\langle u_f(k\delta),\psi_0\rangle_m-\langle u_f((k+1)\delta),\psi_0\rangle_m}
{\langle u_f((k+1)\delta),\psi_0\rangle_m\langle u_f(k\delta),\psi_0\rangle_m}= \delta\nu
$$
uniformly in $f\in\C_1^+\setminus\{0\}$.
Now, \eqref{2.3} follows immediately.
\hfill$\Box$

\bigskip

\textbf{Proof of Theorem \ref{lem:2.1}:}
For $t>0$, we have
\begin{equation}\label{4.1}
{\P}_{\mu}\left( \|X_t\|\ne 0 \right)=
 \lim_{\theta\to\infty}\left(1-\exp\{-\langle u_{\theta}(t),\mu\rangle\}\right).
\end{equation}
Using Lemma \ref{le:1.4} with $\delta=1$, we have
\begin{equation}\label{2.13}
  \lim_{n\to\infty}\frac{1}{n}\left(\frac{1}{\langle u_{\theta}(n),\psi_0\rangle_m}-\frac{1}{\theta\langle1 ,\psi_0\rangle_m}\right)=\nu
\end{equation}
uniformly in $\theta>0$.
For $\theta>1$, it holds that
\begin{equation}\label{2.23}
  \frac{1}{n}\frac{1}{\theta\langle 1,\psi_0\rangle_m}\le \frac{1}{n}\frac{1}{\langle 1,\psi_0\rangle_m}\to 0,
  \quad \mbox{ as } n\to\infty,
\end{equation}
uniformly in $\theta>1$.
It follows from \eqref{2.13} and \eqref{2.23} that
\begin{equation}\label{2.24}
  \lim_{n\to\infty}n\langle u_{\theta}(n),\psi_0\rangle_m=\nu^{-1}
\end{equation}
uniformly in $\theta>1$.
By \eqref{1.10} and \eqref{1.22}, we have, as $n\to\infty$,
for any $\mu(E)\le M$,
\begin{eqnarray*}
  &&n|\langle  u_{\theta}(n),\mu\rangle-\langle  u_{\theta}(n),\psi_0\rangle_m\langle \phi_0,\mu\rangle|= n\langle  u_{\theta}(n),\psi_0\rangle_m|\langle h_\theta(n)\phi_0,\mu\rangle|\nonumber\\
  &\le& M\|h_{\theta}(n)\|_\infty \phi_0\|_\infty n\langle  u_{\theta}(n),\psi_0\rangle_m\|\to 0,
\end{eqnarray*}
uniformly in $\theta>1$.
Thus,
\begin{equation}\label{4.13}
  \lim_{n\to\infty}n \langle  u_{\theta}(n),\mu\rangle=\nu^{-1}\langle \phi_0,\mu\rangle\quad
  \mbox{uniformly in }\theta>1 \mbox{ and }\mu \mbox{ with } \mu(E)\le M.
\end{equation}
By \eqref{2.5}, we have $\langle  u_{\theta}(n),\mu\rangle\le \langle -\log q_{n},\mu\rangle\le \|-\log q_{n}\|_\infty\|\mu\|\to 0$, as $n\to\infty$, uniformly
in $\theta>0$ and $\mu$ with $\mu(E)\le M$.
Therefore, it follows from \eqref{4.13} that
$$\lim_{n\to\infty}n\left(1-\exp\{-\langle  u_{\theta}(n),\mu\rangle\}\right)=\nu^{-1}\langle\phi_0,\mu\rangle
 \mbox{ uniformly in }\theta>1 \mbox{ and }\mu \mbox{ with } \mu(E)\le M.
$$
Hence by \eqref{4.1}, we have
$$
  \lim_{n\to\infty}n{\P}_{\mu}\left( \|X_{n}\|\ne 0 \right)=\nu^{-1}\langle\phi_0,\mu\rangle,
$$
uniformly in $\mu$ with $\mu(E)\le M$.
Since ${\P}_{\mu}\left( \|X_{t}\|\ne 0 \right)$ is decreasing in $t$, we have
$$
[t] \, {\P}_{\mu}\left( \|X_{([t]+1)}\|\ne 0 \right)\le t\,{\P}_{\mu}\left( \|X_{t}\|\ne 0 \right)
\le ([t]+1)\, {\P}_{\mu}\left( \|X_{[t]}\|\ne 0 \right).
$$
Now \eqref{2.4} follows immediately.
\hfill$\Box$

\bigskip
Now we are ready to prove  Theorem \ref{The:2.2}.

\textbf{Proof of Theorem \ref{The:2.2}:}
First, we consider the special case when $f(x)=\phi_0(x)$.
We only need to show that, for any $\lambda>0$,
\begin{equation}\label{2.33}
  \P_{\mu}\left(\exp\left\{-\lambda t^{-1}\langle\phi_0,X_t\rangle\right\}\mid \|X_t\|\ne 0\right)
  \to\frac{1}{\lambda\nu+1}, \quad \mbox{ as } t\to\infty.
\end{equation}
Note that
\begin{eqnarray*}
  &&\P_{\mu}\left(\exp\left\{-\lambda t^{-1}\langle\phi_0,X_t\rangle\right\}\mid \|X_t\|\ne 0\right)\\
  &&=\frac{\P_{\mu}\left(\exp\left\{-\lambda t^{-1}\langle\phi_0,X_t\rangle\right\}\right)-\P_{\mu}(\|X_t\|= 0)}{\P_{\mu}(\|X_t\|\ne 0)}\\
  &&= 1-\frac{1-\P_{\mu}\left(\exp\left\{-\lambda t^{-1}\langle\phi_0,X_t\rangle\right\}\right)}{\P_{\mu}(\|X_t\|\ne 0)}.
  \end{eqnarray*}
By Theorem \ref{lem:2.1}, to prove \eqref{2.33}, it suffices to show that, as $t\to\infty$,
\begin{equation}\label{2.31}
  t\left(1-\P_{\mu}\left(\exp\left\{-\lambda t^{-1}\langle\phi_0,X_t\rangle\right\}\right)\right)
  =t \left(1-\exp\left\{-\langle  u_{\lambda t^{-1}\phi_0}(t),\mu\rangle\right\}\right)\to \frac{\lambda}{\lambda\nu+1}\langle\phi_0,\mu\rangle.
\end{equation}
Since $X_t$ is right continuous and $\phi_0$ is a bounded continuous function, $t\rightarrow \P_{\mu}\left(\exp\left\{-\lambda t^{-1}\langle\phi_0,X_t\rangle\right\}\right)$ is a right continuous function.
By the Croft-Kingman lemma (see, for example, \cite[Section 6.5]{WA}),
it suffices to show that, for every $\delta>0$, \eqref{2.31} holds for every sequence $n\delta$ as $n\to\infty$.
For this, it is enough to prove that for any $\delta>0$, as $n\to\infty$,
  \begin{equation}\label{4.3}
  n\delta\,\langle  u_{\lambda (n\delta)^{-1}\phi_0}(n\delta),\mu\rangle\to \frac{\lambda}{\lambda\nu+1}\langle \phi_0,\mu\rangle.
\end{equation}
By Lemma \ref{le:1.4}, we have
\begin{eqnarray*}
  &&\lim_{n\to\infty}\frac{1}{(n\delta)\langle  u_{\lambda (n\delta)^{-1}\phi_0}(n\delta),\psi_0\rangle_m} \\
  &=& \lim_{n\to\infty}\frac{1}{n\delta}\left(\frac{1}{\langle  u_{\lambda (n\delta)^{-1}\phi_0}(n\delta),\psi_0\rangle_m}-\frac{1}
  {\langle\lambda(n\delta)^{-1}\phi_0,\psi_0\rangle_m}\right)
  +\frac{1}{\lambda} \\
   &=&  \nu+\lambda^{-1},
\end{eqnarray*}
which implies that
\begin{equation}\label{2.32}
  (n\delta)\langle  u_{\lambda (n\delta)^{-1}\phi_0}(n\delta),\psi_0\rangle_m\to \frac{\lambda}{\lambda\nu+1},
  \quad \mbox{ as } n\to\infty.
\end{equation}
Using Lemma \ref{lem:2.3} and \eqref{2.32}, we get that, as $n\to\infty$,
\begin{eqnarray}\label{4.5}
  &&n\delta\left|\langle  u_{\lambda (n\delta)^{-1}\phi_0}(n\delta),\mu\rangle-\langle  u_{\lambda (n\delta)^{-1}\phi_0}(n\delta),\psi_0\rangle_m\langle \phi_0,\mu\rangle\right|\nonumber\\
  &\le& n\delta\|h_{\lambda (n\delta)^{-1}\phi_0}(n\delta)\|_\infty\langle  u_{\lambda (n\delta)^{-1}\phi_0}(n\delta),\psi_0\rangle_m\langle \phi_0,\mu\rangle\to 0.
\end{eqnarray}
Now \eqref{4.3} follows easily from \eqref{2.32} and \eqref{4.5}.

For a general $f$, let
\begin{equation}\label{tildef}
  \tilde{f}(x)=f(x)-\langle f,\psi_0\rangle_m\phi_0(x).
\end{equation}
Then, $\langle \tilde{f},\psi_0\rangle_m=0$.
It is clear that
\begin{eqnarray}\label{4.14}
  &&\P_{\mu}\left(\left(t^{-1}\langle\tilde{f},X_t\rangle\right)^2\mid \|X_t\|\ne 0\right)
  = \frac{\P_{\mu}\left(\langle\tilde{f},X_t\rangle\right)^2}{t^2\P_{\mu}(\|X_t\|\ne 0)}.
\end{eqnarray}
By the branching property and \eqref{2.11}, we have,
$$\sup_{t>2}\V{\rm ar}_\mu \langle\tilde{f},X_t\rangle=\sup_{t>2}\langle \V{\rm ar}_{\delta_{\cdot}} \langle\tilde{f},X_t\rangle,\mu \rangle<\infty. $$
It follows from \eqref{1.6} that
$$\sup_{t>1}\left|\P_{\mu}\langle\tilde{f},X_t\rangle\right|=\sup_{t>1}\left|\langle T_t\tilde{f},\mu\rangle\right|<\infty.$$
Combining the last two displays, we get that $\sup_{t>2}\P_{\mu}\left(\langle\tilde{f},X_t\rangle\right)^2<\infty$.
Thus by \eqref{2.4} and \eqref{4.14}, we get that as $t\to\infty$,
$$
\P_{\mu}\left(\left(t^{-1}\langle\tilde{f},X_t\rangle\right)^2\mid \|X_t\|\ne 0\right)\to 0,
$$
which implies that, for any
 $\epsilon>0$,
\begin{equation}\label{2.51}
\lim_{t\to\infty}\P_{t,\mu}\left(\left|t^{-1}\langle \tilde{f},X_t\rangle \right|\ge \epsilon\right)=0.
\end{equation}
Thus, by \eqref{tildef}, we have
$$t^{-1}\langle f,X_t\rangle|_{\P_{t,\mu}} \stackrel{d}{\rightarrow}\langle f,\psi_0\rangle_mW.$$
\hfill$\Box$

\textbf{Proof of Corollary \ref{cor1.1}:}
Recall that for $f\in\C_2$, $\tilde{f}$ was defined in \eqref{tildef}. Thus
$$\frac{\langle f,X_t\rangle}{\langle \phi_0,X_t\rangle}-\langle f,\psi_0\rangle_m=\frac{\langle \tilde{f},X_t\rangle}{\langle \phi_0,X_t\rangle}.$$
For any $\epsilon>0$ and $\delta>0$, by \eqref{2.51} and \eqref{the:1}, we have,
\begin{eqnarray*}
  &&\P_{\mu}\left(\frac{|\langle \tilde{f},X_t\rangle|}{\langle \phi_0,X_t\rangle}>\epsilon\mid\|X_t\| \ne 0\right) \\
   &&\le \P_{\mu}\left(t^{-1}|\langle \tilde{f},X_t\rangle|>\delta\mid\|X_t\| \ne 0\right)
   +\P_{\mu}\left(t^{-1}\langle \phi_0,X_t\rangle<\delta/\epsilon\mid\|X_t\| \ne 0\right)\\
   &&\to 0+P\left(W<\delta/\epsilon\right),\quad \mbox{as }t\to\infty.
\end{eqnarray*}
Letting $\delta\to0$, we get that
$$\lim_{t\to\infty}\P_{t,\mu}\left(\frac{|\langle \tilde{f},X_t\rangle|}{\langle \phi_0,X_t\rangle}>\epsilon\right)=0,$$
which implies \eqref{3.1}.
All real-valued continuous functions with compact support in $E$ belong to $\C_2$.
Thus, by \eqref{3.1}, we have that for any real-valued continuous function $f$ with compact support,
\begin{equation}\label{3.1'}
 \frac{\langle f,X_t\rangle}{\langle \phi_0,X_t\rangle}|_{\P_{t,\mu}}\stackrel{d}{\rightarrow}l(f)(=\langle f,\psi_0\rangle_m),
\end{equation}
Hence by \cite[Theorem 16.16]{Kall}, we get that
\begin{equation*}
 \frac{X_t}{\langle \phi_0,X_t\rangle}|_{\P_{t,\mu}}\stackrel{d}{\rightarrow}l.
\end{equation*}
Since $\nu\to \rho(\nu, l)\wedge 1$ is a bounded continuous function on the
space of Radon measures on $E$ equipped with the vague topology, we have
$$
\lim_{t\to\infty}\P_{t,\mu}\left[\rho\left(\frac{X_t}{\langle \phi_0,X_t\rangle} ,l\right)\wedge1\right]=0,
$$
from which the last assertion of the corollary follows immediately.
\hfill$\Box$

\subsection{Proof of Theorem \ref{them2}}

In this subsection, we give the proof of Theorem \ref{them2}. We prove a simple lemma first.

\begin{lemma}\label{lem3.2}
Suppose that $\mathcal{V}$ is an index set
and $\{F_v: v\in \mathcal{V}\}$ is a family of uniformly bounded random variables, that is, there is a constant $M$ such that $|F_v|\le M$ for all $v\in\mathcal{V}$, then any $s>0$,
\begin{equation}\label{3.4}
\lim_{t\to\infty}\sup_{v\in\mathcal{V}}|\P_{t+s,\mu}(F_v)-\P_{t,\mu}(F_v)|= 0.
\end{equation}
\end{lemma}
\textbf{Proof:}
By Theorem \ref{lem:2.1}, we have
\begin{equation}\label{3.5}
\lim_{t\to\infty}\frac{\P_{\mu}(\|X_{t}\|\ne 0)}{\P_{\mu}(\|X_{t+s}\|\ne 0)}=1.
\end{equation}
By the definition of $\P_{t,\mu}$, we have
\begin{eqnarray*}
  \P_{t+s,\mu}(F_{v}) &=& \P_{t,\mu}(F_{v}, \ \|X_{t+s}\|\ne0)\frac{\P_{\mu}(\|X_{t}\|\ne 0)}
  {\P_{\mu}(\|X_{t+s}\|\ne 0)} \\
   &=&\P_{t,\mu}(F_{v})\frac{\P_{\mu}(\|X_{t}\|\ne 0)}{\P_{\mu}(\|X_{t+s}\|\ne 0)}
   -\P_{t,\mu}(F_{v},\|X_{t+s}\|=0)\frac{\P_{\mu}(\|X_{t}\|\ne 0)}{\P_{\mu}(\|X_{t+s}\|\ne 0)}.
\end{eqnarray*}
Thus, as $t\to \infty$,
\begin{eqnarray*}
  &&|\P_{t+s,\mu}(F_{v})-\P_{t,\mu}(F_{v})|
  \le M\left|\frac{\P_{\mu}(\|X_{t}\|\ne 0)}
  {\P_{\mu}(\|X_{t+s}\|\ne 0)}-1\right|+M\P_{t,\mu}(\|X_{t+s}\|=0)\frac{\P_{\mu}(\|X_{t}\|\ne 0)}{\P_{\mu}(\|X_{t+s}\|\ne 0)}\\
    &=& M\left|\frac{\P_{\mu}(\|X_{t}\|\ne 0)}
  {\P_{\mu}(\|X_{t+s}\|\ne 0)}-1\right|+M\frac{\P_{\mu}(\|X_{t+s}\|=0,\|X_t\|\ne 0)}{\P_{\mu}(\|X_t\|\ne0)}\frac{\P_{\mu}(\|X_{t}\|\ne 0)}{\P_{\mu}(\|X_{t+s}\|\ne 0)}\\
  &=& M\left|\frac{\P_{\mu}(\|X_{t}\|\ne 0)}
  {\P_{\mu}(\|X_{t+s}\|\ne 0)}-1\right|+M\frac{\P_{\mu}(\|X_t\|\ne 0)-\P_{\mu}(\|X_{t+s}\|\ne 0)}{\P_{\mu}(\|X_{t+s}\|\ne0)}\\
     &=&  2M\left|\frac{\P_{\mu}(\|X_{t}\|\ne 0)}{\P_{\mu}(\|X_{t+s}\|\ne 0)}-1\right|\to 0.
\end{eqnarray*}
\hfill$\Box$

We now recall some facts about weak convergence which will be used later.
For $f:\mathbb{R}^d\to\mathbb{R}$, let $\|f\|_L:=\sup_{x\ne y}|f(x)-f(y)|/\|x-y\|$ and
  $\|f\|_{BL}:=\|f\|_{\infty}+\|f\|_L$.
   For any probability measures $\nu_1$ and $\nu_2$ on $\mathbb{R}^d$, define
\begin{equation*}
  \beta(\nu_1,\nu_2):=\sup\left\{\left|\int f\,d\nu_1-\int f\,d\nu_2\right|~:~\|f\|_{BL}\leq1\right\}.
\end{equation*}
Then $\beta$ is a metric.
It follows from \cite[Theorem 11.3.3]{Dudley} that the topology generated by $\beta$
is equivalent to the weak convergence topology.
 From the definition, we can easily see that, if $\nu_1$ and $\nu_2$ are the distributions of two $\R^d$-valued random variables $X$ and $Y$ respectively, then
\begin{equation}\label{5.20}
  \beta(\nu_1,\nu_2)\leq E\|X-Y\|\leq\sqrt{ E\|X-Y\|^2}.
\end{equation}

The following simple fact will be used several times later in this section:
\begin{equation}\label{3.20}
  \left|e^{ix}-\sum_{m=0}^n\frac{(ix)^m}{m!}\right|\leq \min\left(\frac{|x|^{n+1}}{(n+1)!}, \frac{2|x|^n}{n!}\right).
\end{equation}

\bigskip

Now we are ready to prove Theorem \ref{them2}.

\textbf{Proof of Theorem \ref{them2}:}
Define an $\R^2$-valued random variable:
$$
U_1(t):=\left(t^{-1}\langle \phi_0, X_{t}\rangle,t^{-1/2}\langle f,X_{t}\rangle\right).
$$
We need to prove that,  conditioning on $\|X_t\|\ne 0$, $U_1(t)$ converges
to $\left(W,G(f)\sqrt{W}\right)$ in distribution as $t\to\infty$,
which is equivalent to proving that, when one lets $t$ tend to $\infty$ first
and then lets $s$ tend to $\infty$,
\begin{equation}\label{to prove U1}
U_1(t+s)|_{\P_{t+s,\mu}}\stackrel{d}{\to}\left(W,G(f)\sqrt{W}\right).
\end{equation}

Before we prove \eqref{to prove U1}, we first give the main idea of the proof.
In Theorem \ref{The:2.2}, we have proved that the first component of $U_1(t)$ converges to $W$.
So the key is the second component. If we condition on $X_t$,
the mean of $\langle f,X_{s+t}\rangle$ is $\langle T_sf,X_t\rangle$.
Let us consider the centered random variable
$\langle f,X_{s+t}\rangle-\langle T_sf,X_t\rangle$. For fixed $s>0$,
as $t\to\infty$, since the `infinitesimal particles' evolve
independently after time $t$, it is reasonable to expect that,
conditioning on $X_t$ and $\|X_t\|\neq0$,
$\left(\langle f,X_{s+t}\rangle-\langle T_sf,X_t\rangle\right)/\sqrt{\left(\V{\rm ar}
\langle f,X_{s+t}\rangle|X_t\right)}$ converges in distribution to a standard
normal random variable. Note that
$\V{\rm ar}(\langle f,X_{s+t}\rangle|X_t)=\langle \V{\rm ar}_{\delta_\cdot}
\langle f,X_s\rangle, X_t\rangle$. By Theorem \ref{The:2.2}, we have
$t^{-1}\V{\rm ar}(\langle f,X_{s+t}\rangle|X_t)\stackrel{d}{\to}
\langle \V{\rm ar}_{\delta_\cdot}\langle f,X_s\rangle,
\psi_0\rangle_m W$ as $t\to\infty$. We may thus conclude that
$t^{-1/2}\left(\langle f,X_{s+t}\rangle-\langle T_sf,X_t
\rangle\right)\stackrel{d}{\to}\sqrt{W}G_s$,  where
$G_s\sim\mathcal{N}(0,\sigma^2_f(s))$
 with $\sigma^2_f(s)=\langle
\V{\rm ar}_{\delta_\cdot}\langle f,X_s\rangle,\psi_0\rangle_m$ and $W$
is the random variable defined in Theorem \ref{The:2.2}.

The above analysis suggests that we should first consider another $\R^2$-valued random
variable $U_2(s,t)$ defined by
$$
U_2(s,t)=\left(t^{-1}\langle \phi_0,X_{t}\rangle, t^{-1/2}\left(\langle f,X_{s+t}\rangle-\langle T_sf,X_t\rangle\right)\right)\quad s,t>2.
$$
We claim that,
\begin{equation}\label{10.5}
     U_2(s,t)|_{\P_{t,\mu}}\stackrel{d}{\to}\left(W, \sqrt{W}G_s\right), \quad \mbox{ as }
     t\to\infty.
\end{equation}
We will leave the proof of \eqref{10.5} to the end of the proof of this theorem.

Define
$$
U_3(s,t):=\left((t+s)^{-1}\langle \phi_0,X_{t}\rangle, (t+s)^{-1/2}\left(\langle f,X_{s+t}\rangle-\langle T_sf,X_t\rangle\right)\right).
$$
By \eqref{10.5}, we have
\begin{equation}\label{4.10}
U_3(s,t)|_{\P_{t,\mu}}\stackrel{d}{\to}
   \left(W, \sqrt{W}G_s\right),
\end{equation}
as $t\to\infty$.
It follows from \eqref{T_t} and \eqref{2.4} that, as $t\to\infty$,
\begin{eqnarray*}
  (t+s)^{-2}\P_{t,\mu}\left(\langle \phi_0,X_{t+s}\rangle-\langle \phi_0,X_{t}\rangle\right)^2 &=& \frac{\P_{\mu}\left(\langle \phi_0,X_{t+s}\rangle-\langle \phi_0,X_{t}\rangle\right)^2 }{(t+s)^2\P_{\mu}(\|X_t\|\ne 0)} \\
   &=& \frac{\P_{\mu}\left(\langle \V{\rm ar}_{\delta_{\cdot}}\langle \phi_0,X_s\rangle,X_{t}\rangle\right) }{(t+s)^2\P_{\mu}(\|X_t\|\ne 0)}\to 0.
\end{eqnarray*}
If we put
$$
U_4(s,t):=\left((t+s)^{-1}\langle \phi_0,X_{t+s}\rangle, (t+s)^{-1/2}\left(\langle f,X_{s+t}\rangle-\langle T_sf,X_t\rangle\right)\right),
$$
then, as $t\to\infty$,
\begin{equation}\label{4.11}
  U_4(s,t)|_{\P_{t,\mu}}\stackrel{d}{\to}
   \left(W, \sqrt{W}G_s\right).
\end{equation}
By Lemma \ref{lem3.2}, we have
\begin{eqnarray*}
  &&\lim_{t\to\infty}\P_{t+s,\mu}
  \left(\exp\left\{i\theta_1(t+s)^{-1}\langle \phi_0,X_{t+s}\rangle+i\theta_2(t+s)^{-1/2}\left(\langle f,X_{s+t}\rangle-\langle T_sf,X_t\rangle\right)\right\}\right)\\
  &=& \lim_{t\to\infty}\P_{t,\mu}
  \left(\exp\left\{i\theta_1(t+s)^{-1}\langle \phi_0,X_{t+s}\rangle+i\theta_2(t+s)^{-1/2}\left(\langle f,X_{s+t}\rangle-\langle T_sf,X_t\rangle\right)\right\}\right)\\
  &=&P\left(\exp\left\{i\theta_1W+i\theta_2\sqrt{W}G_s\right\}\right).
\end{eqnarray*}
Thus, we have
\begin{equation}\label{4.12}
  U_4(s,t)|_{\P_{t+s,\mu}}\stackrel{d}{\to}
   \left(W, \sqrt{W}G_s\right), \quad \mbox{ as } t\to\infty.
\end{equation}
Now, we deal with $J_2(t,s):=\frac{\langle T_sf,X_t\rangle}{(t+s)^{1/2}}$.
We claim that
\begin{equation}\label{3.10}
  \lim_{s\to\infty}\limsup_{t\to\infty}\P_{t+s,\delta_x}\left(|J_2(t,s)|^2\right)=0.
\end{equation}
By \eqref{T_t}, we have that $\P_{\mu}\langle T_sf,X_t\rangle=\langle T_{t+s}f,\mu\rangle\to 0$ as $t\to\infty$.
Thus by \eqref{2.4} and \eqref{2.10}, we have
\begin{eqnarray}
  &&\limsup_{t\to\infty}\P_{t+s,\delta_x}\left(|J_2(t,s)|^2\right)=  \limsup_{t\to\infty}\frac{\P_{\mu}\Big(\langle T_sf,X_t\rangle^2,\|X_{t+s}\|\neq0\Big)}
   {(t+s)\P_{\mu}(\|X_{t+s}\|\ne0)}\nonumber\\
   &\le&\limsup_{t\to\infty}\frac{\P_{\mu}\Big(\langle T_sf,X_t\rangle^2\Big)}
   {(t+s)\P_{\mu}(\|X_{t+s}\|\ne0)} =\nu\sigma_{(T_sf)}^2.\label{3.12}
\end{eqnarray}
It follows from \eqref{sigma} that, as $s\to\infty$,
\begin{eqnarray*}
  \sigma_{(T_sf)}^2&=&\int_s^\infty\langle A(T_u f)^2,\psi_0\rangle_m \,du \to 0.
\end{eqnarray*}
Now \eqref{3.10} follows immediately.

By \eqref{2.10}, we have
$\lim_{s\to\infty}\V ar_{\delta_x}\langle f,X_s\rangle= \sigma^2_f\phi_1(x)$,
thus $\lim_{s\to\infty}\sigma_f^2(s)=\sigma^2_f$. Hence,
\begin{equation}\label{10.7}
\lim_{s\to\infty}\beta(G_s,G(f))=0.
\end{equation}
Let $\mathcal{D}(s+t)$ and $\widetilde{\mathcal{D}}(s,t)$ be the distributions of $U_1(s+t)$ and $U_4(s,t)$ under $\P_{t+s,\mu}$ respectively,
and let $\widehat{\mathcal{D}}(s)$ and $\mathcal{D}$
be the distributions of $(W, \sqrt{W}G_s)$
and $(W, \sqrt{W}G(f))$ respectively.
Then, using \eqref{5.20}, we have
\begin{eqnarray}\label{10.11}
  \limsup_{t\to\infty}\beta(\mathcal{D}(s+t),\mathcal{D})&\leq&
  \limsup_{t\to\infty}[\beta(\mathcal{D}(s+t),\widetilde{\mathcal{D}}(s,t))+\beta(\widetilde{\mathcal{D}}(s,t),\widehat{\mathcal{D}}(s))
  +\beta(\widehat{\mathcal{D}}(s),\mathcal{D})]\nonumber\\
 &\leq &\limsup_{t\to\infty}(\sqrt{\mathbb{P}_{t+s,\mu}((t+s)^{-1}\langle T_s f,X_{t}\rangle^2)}+0+\beta(\widehat{\mathcal{D}}(s),\mathcal{D}).
\end{eqnarray}
Then we have
$$ \limsup_{t\to\infty}\beta(\mathcal{D}(t),\mathcal{D})= \limsup_{t\to\infty}\beta(\mathcal{D}(s+t),\mathcal{D})\le
\limsup_{t\to\infty}(\sqrt{\mathbb{P}_{t+s,\mu}(J_2(s,t)^2)}+\beta(\widehat{\mathcal{D}}(s),\mathcal{D}).$$
Letting $s\to\infty$, by \eqref{3.10} and \eqref{10.7}, we get
$$\limsup_{t\to\infty}\beta(\mathcal{D}(t),\mathcal{D})=0,$$
which implies the result of theorem.

Now we prove \eqref{10.5}.

Denote the characteristic function of $U_2(s,t)$ under $\mathbb{P}_{t,\mu}$ by
$\kappa_1(\theta_1,\theta_2,s,t)$:
\begin{eqnarray}\label{kappa}
  &&\kappa_1(\theta_1,\theta_2,s,t)\nonumber\\
  &=&\P_{t,\mu}(\exp\{i\theta_1t^{-1}\langle \phi_0,X_{t}\rangle+i\theta_2t^{-1/2}\left(\langle f,X_{s+t}\rangle-\langle T_sf,X_t\rangle\right)\}\nonumber \\
  &=& \P_{t,\mu}\left(\exp\left\{i\theta_1t^{-1}\langle \phi_0,X_{t}\rangle\right.\right.\nonumber\\
  &&\left.\left.+\int_E\int_\mathbb{D}\left( e^{i\theta_2t^{-1/2}\langle f,\omega_s\rangle}-1-i\theta_2t^{-1/2}\langle f,\omega_s\rangle\right)\N_x(d\omega)X_t(dx)\right\}\right),
\end{eqnarray}
where in the last equality we used the Markov property of $X$, \eqref{cf} and \eqref{N1}.
Define
$$J_s(\theta,x):=\int_{\mathbb{D}}\left(\exp\{\langle i\theta f,\omega_s\rangle\}-1-i\theta\langle f,\omega_s\rangle\right)\mathbb{N}_x(d\omega)$$
and
$$I_s(\theta,x):=\int_{\mathbb{D}}\left(\exp\{\langle i\theta f,\omega_s\rangle\}-1-i\theta\langle f,\omega_s\rangle+\frac{1}{2}\theta^2\langle f,\omega_s\rangle^2\right)\mathbb{N}_x(d\omega).$$
Let $V_s(x)=\V{\rm ar}_{\delta_x}\langle f,X_s\rangle\in\C_2^+$. Then, by \eqref{N2}, we have
\begin{eqnarray*}
  J_s(\theta,x) &=& -\frac{1}{2}\theta^2V_s(x)+I_s(\theta,x) \\
   &=& -\frac{1}{2}\theta^2\langle V_s,\psi_0\rangle_m\phi_0(x)-\frac{1}{2}\theta^2\widetilde{V_s}(x)+I_s(\theta,x),
\end{eqnarray*}
where $\widetilde{V}_s=V_s-\langle V_s,\psi_0\rangle_m\phi_0(x)\in\C_2$.
Thus, we have
\begin{eqnarray}\label{4.8}
  &&i\theta_1t^{-1}\langle \phi_0,X_t\rangle+\langle J_s(t^{-1/2}\theta_2,\cdot),X_t\rangle\nonumber \\
   &=& \left(i\theta_1-\frac{1}{2}\theta^2_2\langle V_s,\psi_0\rangle_m\right)t^{-1}\langle \phi_0,X_t\rangle-\frac{1}{2}\theta_2^2t^{-1}\langle \widetilde{V}_s,X_t\rangle+\langle I_s(t^{-1/2}\theta_2,\cdot),X_t\rangle.
\end{eqnarray}
By \eqref{2.51}, we know that, for any $\epsilon>0$,
\begin{equation}\label{4.6}
\lim_{t\to\infty}\P_{t,\mu}\left(\left|t^{-1}\langle \widetilde{V}_s,X_t\rangle\right|\ge \epsilon\right)
=0.
\end{equation}
By \eqref{3.20}, we have
\begin{eqnarray}\label{4.7}
\left|I_s(t^{-1/2}\theta_2,x)\right|
&\leq& \theta_2^2t^{-1}\mathbb{N}_{x}
\left(\langle f,\omega_s\rangle^2\left(\frac{t^{-1/2}\theta_2\langle f,\omega_s\rangle}{6}\wedge 1\right)\right).
\end{eqnarray}
Let
$$h(x,s,t)=\mathbb{N}_{x}
\left(\langle f,\omega_s\rangle^2\left(\frac{t^{-1/2}\theta_2\langle f,\omega_s\rangle}{6}\wedge 1\right)\right).$$
We note that $h(x,s,t)\downarrow0$ as  $t\uparrow\infty$.
By \eqref{2.11}, we have
$$
h(x,s,t)\le \mathbb{N}_{x}(\langle f,X_s\rangle^2)=\mathbb{V}{\rm ar}_{\delta_x}\langle f,X_s\rangle\lesssim \phi_0(x)\in \C_2.
$$
Thus, by \eqref{2.4} and \eqref{T_t}, we have, for any $u<t$,
$$t^{-1}\P_{t,\mu}\langle h(\cdot,s,t),X_t\rangle\le t^{-1}\P_{t,\mu}\langle h(\cdot,s,u),X_t\rangle=\frac{\P_{\mu}\langle h(\cdot,s,u),X_t\rangle}{t\P_{\mu}(\|X_t\|\ne 0)}\to \nu\langle h(\cdot,s,u),\psi_0\rangle_m,$$
as $t\to\infty$.
Letting $u\to\infty$, we get $\langle h(\cdot,s,u),\psi_0\rangle_m\to 0.$
Thus, by \eqref{4.7}, we get that
$$\lim_{t\to\infty}\P_{t,\mu}|\langle I_s(t^{-1/2}\theta_2,\cdot),X_t\rangle|=0,$$
which implies that, for any $\epsilon>0$,
\begin{equation}\label{4.9}
\lim_{t\to\infty} \P_{t,\mu}\left(\left|\langle I_s(t^{-1/2}\theta_2,\cdot),X_t\rangle\right|\ge \epsilon\right)
=0.
\end{equation}
Thus, by \eqref{4.6}, \eqref{4.9} and \eqref{4.8}, we get
$$i\theta_1t^{-1}\langle \phi_0,X_t\rangle+\langle J_s(t^{-1/2}\theta_2,\cdot),X_t\rangle|_{\P_{t,\mu}}\stackrel{d}{\to}\left(i\theta_1
  -\frac{1}{2}\theta_2^2\langle V_s,\psi_0\rangle_m\right)W.$$
Since the real part of $J_s(t^{-1/2}\theta_2,x)$ is non-positive, we have
$$|\exp\{i\theta_1t^{-1}\langle \phi_0,X_t\rangle+\langle J_s(t^{-1/2}\theta_2,\cdot),X_t\rangle\}|\le1.$$
Therefore, by \eqref{kappa} and the dominated convergence theorem, we get
$$\lim_{t\to\infty}\kappa_1(\theta_1,\theta_2,s,t)=P\left(\exp\left\{\left(i\theta_1
  -\frac{1}{2}\theta_2^2\langle V_s,\psi_0\rangle_m\right)W\right\}\right),$$
which implies our claim \eqref{10.5}.
\hfill$\Box$

\smallskip

\smallskip

\begin{singlespace}

\end{singlespace}

\vskip 0.2truein
\vskip 0.2truein

\noindent{\bf Yan-Xia Ren:} LMAM School of Mathematical Sciences \& Center for
Statistical Science, Peking
University,  Beijing, 100871, P.R. China. Email: {\texttt
yxren@math.pku.edu.cn}

\smallskip
\noindent {\bf Renming Song:} Department of Mathematics,
University of Illinois,
Urbana, IL 61801, U.S.A.
Email: {\texttt rsong@math.uiuc.edu}

\smallskip

\noindent{\bf Rui Zhang:} School of Mathematical Sciences, Capital Normal
University,  Beijing, 100048, P.R. China. Email: {\texttt
zhangrui27@cnu.edu.cn}

\end{doublespace}

\end{document}